\documentclass[12pt]{article}
\usepackage[margin=0.7in]{geometry}
\usepackage{amssymb,amsmath,amsthm,mathrsfs}
\usepackage{graphicx,epstopdf,color}
\usepackage{array}

\allowdisplaybreaks
\usepackage[alphabetic]{amsrefs}
\usepackage [latin1]{inputenc}
\newcommand{\R}{\mathbb{R}}
\newcommand{\Z}{\mathbb{Z}}

\newcommand{\C}{\mathbb{C}}

\renewcommand{\leq}{\leqslant}
\renewcommand{\le}{\leqslant}
\renewcommand{\geq}{\geqslant}
\renewcommand{\ge}{\geqslant}

\usepackage{hyperref}
\usepackage{authblk}

\begin{document}

\title{Efficiency functionals for
the L\'evy flight foraging hypothesis}

\author{Serena Dipierro, Giovanni Giacomin and Enrico Valdinoci
\thanks{Department of Mathematics and Statistics,
University of Western Australia, 35 Stirling Highway,
WA6009 Crawley, Australia.\\
{\tt serena.dipierro@uwa.edu.au,
giovanni.giacomin.1@studenti.unipd.it,
enrico.valdinoci@uwa.edu.au}}}

\maketitle

\begin{abstract}
We consider a forager diffusing via a fractional
heat equation and we introduce several efficiency functionals
whose optimality is discussed in relation to the
L\'evy exponent of the evolution equation.

Several biological scenarios, such as
a target close to the forager, a sparse environment,
a target located away from the forager and two targets
are specifically taken into account.

The optimal strategies of each of these
configurations are here analyzed explicitly
also with the aid of some special functions
of classical flavor and the results are confronted
with the existing paradigms of the
L\'evy foraging hypothesis.

Interestingly, one discovers bifurcation phenomena
in which a sudden switch occurs between
an optimal (but somehow unreliable) L\'evy foraging pattern
of inverse square law
type and a less ideal (but somehow more secure) classical Brownian motion
strategy. 

Additionally, optimal foraging strategies can be detected in
the vicinity of the Brownian one even in cases in which the
Brownian one is pessimizing an efficiency functional.
\end{abstract}

\section{Introduction}

Foraging theory (see e.g.~\cite{STEP}) is a fascinating, important and cross-disciplinary
topic of investigation that gathers
together researchers from different areas (such as biologists, ethologists, physicists,
statisticians, computer scientists, mathematicians, etc.). It is commonly accepted that
the broad variety of environmental and biological situations in nature and the
Darwinistic
evolution through natural selection have led over time to highly
efficient foraging strategies
see e.g.~\cite{QUANTI}
(it is however under an intense debate whether L\'evy type
patterns in animal searches are an evolutionary stable and
well consolidated outcome~\cite{H1}
or they are produced by innate composite correlated random walks~\cite{MUSS0, MUSS1};
under investigation
is also the role of particular distribution of resources for the emergence of foraging patterns, see e.g.~\cite{PRIM};
it is also
debatable that natural selection alone can always
optimize a specific parameter in complex environments, see e.g.~\cite{DARCR, WASTE}).

In general, the precise determination of optimal foraging strategies depends in
a very complex way on a large number of parameters (such as the density and mobility of the preys
and of the searchers and the mutable environmental conditions); furthermore, the collection and analysis of
empirical data are typically challenging tasks, also open to controversial interpretations
due to the use of different mathematical models or even due to spurious information
(see e.g.~\cite{SPU}).

A rather consolidated attempt to understand and classify different foraging strategies according to
the evolution of the distribution of the searchers lies in the so-called L\'evy flight foraging hypothesis.
Namely, rather than diffusing in analogy to the classical Brownian motion, empirical evidence has often
backed the hypothesis that animals move according to a scale-free fractal-like pattern similar to the one
produced by long-jump random walks of L\'evy type,
possibly to avoid being trapped in a
search of food confined in a narrow region beyond
sensory range and to reduce the chances of intensively revisiting immediate
surrounding areas in environments of scarce resources (see e.g.~\cite{SCA, HUM, H1, H2} for empirical
evidence for such biological L\'evy flights).
Phenomena related to L\'evy flights are attracting increasing interest and they seem to possess some
kind of universality, occurring also in situations different from animal foraging and including,
among the others, human settlements and travels, see~\cite{BERT, BCC, BRO, GO, 2014, 2018}
and also~\cite{VIRTU} for related virtual reconstructions.
L\'evy patterns also emerge in dynamical models
as a non-Gaussian transport related to chaos, see e.g.~\cite{CHAOS, STRANGE}.

Several studies have exploited tools from
mathematical analysis and statistical mechanics
to validate the hypothesis that L\'evy flights confer a significant advantage for foragers,
see~\cite{VIS, BART, VIS2, RAP}.
Typically, to confirm the L\'evy flights optimality,
structural assumptions on the environment,
on the searcher and on the target are taken, such as:
the foraging should be of non-destructive type
(that is, once a target has been foraged, it has to reappear infinitely fast);
after foraging,
the seeker starts a new flight ``infinitely close'' to the previous target;
the searcher moves
rapidly relative to the target; the target density is low;
the forager does not keep memory of
previous encounters; the
forager has inadequate information on the area to patrol and
on the target location,
etc.
Of course,
all these characteristics provide a
highly simplified representations of real foragers,
yet conceptual simplifications (rather than trivializations)
are often very advantageous to advance and consolidate
the knowledge on a complex topic.
As a matter of fact,
due to the difficulty of the analytical setting (and also to mimic situations of biological interest),
to develop a mathematical theory of foraging related to the L\'evy flight hypothesis it is often
necessary to introduce additional parameters (such as a ``direct vision distance'' of the predator,
see page~912 in~\cite{VIS}) and approximations (see e.g. equations~(2) and~(5) in~\cite{VIS}).
In general, in spite of several quite strong and convincing attempts to completely deduce the L\'evy flight
foraging hypothesis from prime principles, several important details have generated
debate,
see e.g.~\cite{Palyulin2931, LIT1, LIT2, LIT3}
and also~\cite{JIM} for a review of several
controversial aspects of the L\'evy foraging hypothesis.\medskip

In this paper, we consider the optimal foraging strategies
in several situations of biological interest, such as:
\begin{itemize}
\item the case in which
a single target is located in the proximity of the forager's burrow,
\item the case in which targets are sparsely distributed,
\item the case in which a single
target is located far away from 
the forager's burrow, 
\item the case in which there are two targets,
one close and one far from
the forager's burrow.
\end{itemize}
The optimal strategies of each of these
configurations will be analyzed
in light of new efficiency functionals
relying also on methods from mathematical analysis
and with the aid of some classical special functions.

A few comments are in order to highlight some of the main
structural differences between our approach and the rather abundant
existing literature on optimal animal foraging.
On the one hand, the models considered here
share with the existing literature several common treats, such as the assumption that the forager has no memory about the targets previously hit and that the pray has no awareness of the strategy and the movement of the predator. On the other hand, our models present significant differences with the existing literature for at least the following features:
\begin{itemize}
\item the forager does not restart its strategy after hitting each single target (instead, the seeker diffuses according to a L\'evy type of diffusive equation, and this feature happens to be consistent
with the setting of some of the existing literature, see equation~(1) in~\cite{Palyulin2931};
similar, but different, space-fractional equations in biological environment
have also been considered in view of the Caputo derivative, see equation~(2.1) in~\cite{2121}),
\item no additional parameter related to direct vision is taken into account,
no a-priori bound on step lengths is imposed,
no truncation of the power law distribution is assumed
(with the advantage of not endowing the problem with
auxiliary and sometimes arbitrary parameters;
as a counterpart of these technical and conceptual simplifications,
the diffusions corresponding to infinite mean displacements are ruled
out as infinite overshooting and this feature happens to be consistent
with the setting of some of the existing literature, see e.g.the discussion after formula~(2)
in~\cite{Palyulin2931}),
\item we will take into account time averages of foraging success (though some integrals
over time were previously considered, such as in the cumulative probability in equation~(6) of~\cite{Palyulin2931}, we will specialize our analysis in detecting
different optimal strategies according to the different time scales
involved in the seeking process, rather than simply considering
the foraging outcome at a given time),
\item we will analyze in detail the role played by possibly different normalizing
constants appearing when linking probabilistic models to analytical ones
(typically, these constants depend\footnote{In terms of optimization
strategies, we think it would have been beneficial, for instance, to discuss more
extensively the possible dependences on~$\mu $ (corresponding to~$1+2s$ here)
in the right hand sides of equation~(5) in~\cite{VIS} and equation~(11)
in~\cite{RAP}, as well as the constant~$C$ in equation~(4)
of~\cite{VIUJNA90ok-lm}. The explanation for the pseudo mean squared displacements in the right hand side of equation~(10) of~\cite{RAP} could have also benefitted from further details on the possible
dependence of~$\delta$ and~$\alpha$ (the latter corresponding to~$2s$ here).}
on the fractional exponent~$s$, hence they may play a significant role when the objective is to optimize in~$s$
and, in general, they cannot be
light-heartedly disregarded),
\item we will find solutions in closed form, relying only on elementary special functions (and, since these functions, such as the Euler Gamma Function and the Riemann Zeta Function, are widely studied and already carefully implemented in all mathematical softwares, in our approach no expensive or advanced numerical simulations are needed),
\item we introduce a number of new efficiency functionals
whose optimization can be explicitly discussed (these functionals
are inspired by, but somewhat different from, the mean first passage time adopted in~\cite{VIS} -- in this
way, we also avoid any overlap with
some controversial details in the contemporary literature such
as in~\cite{LIT1, LIT2, LIT3}).
\end{itemize}

Though the arguments developed here essentially
carry over to the multi-dimensional case, for the sake of
simplicity (and following a consolidated tradition in mathematical biology,
see e.g.~\cite{BUL, Palyulin2931}), we stick here to dimension~$1$.
The multi-dimensional case will be treated separately in a subsequent
work, also taking into account new sets of structural parameters
according to the geometry of the space and of the diffusive process.

Also, we focus here on the case of stationary
targets (the case of mobile preys possibly with different velocities
will be accounted for in a forthcoming work).

\begin{center}
\begin{tabular}{ ||p{7.5cm}|p{7.7cm}|| }
 \hline
 \multicolumn{2}{|c|}{{\em Notation Table}} \\
 \hline\hline
 Fourier Transform of $f$ &
 $\widehat{f}(\xi):={\mathcal{F}}f(\xi):=\displaystyle\int_{\R }f(x)e^{-2\pi ix \xi}\,dx$\\
 \hline
 Fourier Antitransform of $g$ &
 $\check{g}(x):={\mathcal{F}}^{-1}g(x):=\displaystyle\int_{\R }g(\xi)e^{2\pi ix \xi}\,d\xi$\\
\hline
Poisson Summation Formula &
 $\displaystyle\sum_{k\in\Z }f(x+k)=\displaystyle\sum_{k\in\Z }\widehat f(k)\,e^{2\pi ix k}$\\
\hline
Dirac Delta Function at~$x_0$ &$ \delta_{x_0}$\\
\hline
Integral of a continuous function~$\phi$ against the Dirac Delta &$\displaystyle
\int_\R \phi(x)\delta_{x_0}(x)\,dx:=\phi(x_0)$\\
\hline
Fractional parameter &$ s\in(0,1)$\\
\hline
Fractional Laplacian of~$u$ &$ (-\Delta)^su:= {\mathcal{F}}^{-1} (|2\pi \xi|^{2s} \widehat u)$\\
\hline
Gamma Function ($z\in\C$, $\Re (z)>0$)
&$\Gamma (z):=\displaystyle\int _{0}^{+\infty }\vartheta^{z-1}e^{-\vartheta}
\,d\vartheta$
\\ \hline
Euler-Mascheroni constant
& $\gamma :=\displaystyle
\lim _{n\to \infty }\left(-\ln n+\sum _{k=1}^{n}{\frac {1}{k}}\right)
=0.5772156...$
\\
\hline
Digamma Function &
$ \displaystyle \psi (z):={\frac {d}{dz}}\ln {\big (}\Gamma (z){\big )}={\frac {\Gamma '(z)}{\Gamma (z)}}
$\\
\hline
Riemann Zeta Function ($z\in\C$, $\Re (z)>1$) &
$\zeta (z):=\displaystyle\sum _{j=1}^{\infty }{\frac {1}{j^{z}}}$
\\
\hline\end{tabular}\end{center}

The results obtained will detect the optimal exponent~$s$
corresponding to the most efficient foraging strategy (according to
the different possible efficiency functionals). Several interesting
patterns will arise. Quite often, optimal strategies are obtained in nature
either by L\'evy flights modeled on the inverse square law,
or by the classical Brownian motion, or by some intermediate fractional values.
In our discussion, all these three patterns will clearly arise
and suitable bifurcation of optimal strategies will occur
in dependence of the environmental parameters.

For instance, varying the time in which the search occurs or the sparseness
of the targets, the optimality of the inverse square law may be lost
in favor of a classical Gaussian strategy (or viceversa), and in some
cases optimal values are found arbitrarily close to pessimal ones
(and, conversely, pessimal values arbitrarily close to optimal ones).
We think that this is a very interesting phenomenon, underlying the fact
that the theoretical optimality of the strategy by itself might be not the main
information to take into account for efficient search algorithms, since
less ideal strategies might produce more consistent results
and prove themselves to be more reliable and viable in concrete situations.

In some circumstances, we will also detect optimal fractional values
of intermediate type between the inverse square law and the Gaussian.
In all cases, we will develop
explicit (and somewhat ``elegant'') representations
of the efficiency functional that we introduce, thus allowing
simple and effective analytic manipulations.
As a byproduct, many of the environmental bifurcation parameters
will be computed exactly.\medskip

The rest of the paper is organized as follows.
In the forthcoming Section~\ref{KS-lSSETC},
we introduce our mathematical setting adopted in this paper,
modeled on a forager randomly diffusing through a nonlocal heat equation
and immobile targets with different types of distributions
(see also the Notation Table for the list
of the main mathematical objects and notations utilized in this paper).
The different biological scenarios corresponding to
these distributions of resources will be discussed
in Section~\ref{DIOFFE-R},
where several efficient functionals will be optimized
with respect to the diffusion exponent. The results
obtained will be also compared with the existing literature
related to the L\'evy foraging hypothesis.

\section{Mathematical setting}\label{KS-lSSETC}

We introduce now the formal mathematical that
setting we work with.
The setting is modeled on the fractional heat equation
and goes as follows.

Let~$\kappa>0$, $s\in(0,1)$ and~$u(x,t)$ be the solution\footnote{In several
occurrences in the existing literature, the L\'evy exponent in biological contexts is denoted by~$\mu$.
With respect to our notation, it holds that~$\mu=1+2s$.}
of
\begin{equation}\label{BC:EQ} \begin{cases}
\partial_t u=-\kappa^{2s}\,(-\Delta)^s u& {\mbox{ in }}\R \times(0,+\infty),\\
u(x,0)=\delta_0(x).
\end{cases}\end{equation}
By taking the Fourier Transform of this relation,
$$ \begin{cases}
\partial_t \widehat u=-|2\pi \kappa\xi|^{2s} \widehat u& {\mbox{ in }}\R \times(0,+\infty),\\
\widehat u(x,0)=1.
\end{cases}$$
Therefore
\begin{equation}\label{WIDE} \widehat u(\xi,t)=\exp\big(-|2\pi \kappa\xi|^{2s}t\big)\qquad{\mbox{and}}\qquad
u(x,t)={\mathcal{F}}^{-1}\Big(\exp\big(-|2\pi \kappa\xi|^{2s}t\big)\Big).
\end{equation}
It is possible that the similarity (and the difference)
between the expression for~$ \widehat u$ in~\eqref{WIDE}
and the standard Gaussian (corresponding to~$s=1$) were one of the inspiring motivations
for L\'evy's approach to the Central Limit Theorem in presence of infinite moments,
see equation~(7) in~\cite{STRANGE}.

We observe that, by scaling,
\begin{equation}\label{SCALE}
\begin{split}
u(x,t)\,&=\,
\int_{\R }
\exp\big(-|2\pi \kappa\xi|^{2s}t+2\pi ix \xi\big)\,d\xi\\&=\,\frac{1}{t^{ \frac{1}{2s}}}
\int_{\R }
\exp\big(-|2\pi\kappa\eta|^{2s} +2\pi i t^{-\frac1{2s}}x \eta\big)\,d\eta\\&=\,\frac{1}{t^{ \frac{1}{2s}}}
{\mathcal{F}}^{-1}\Big(\exp\big(-|2\pi \kappa\xi|^{2s}\big)\Big)\left(\frac{x}{t^{\frac1{2s}}}\right)\\&=\,\frac{1}{t^{ \frac{1}{2s}}}
u\left( \frac{x}{t^{\frac1{2s}}},1\right).
\end{split}
\end{equation}
In addition (see e.g. formula~(2.30) in~\cite{ABA}),
\begin{equation}\label{GIU}
0\le u(x,1)\le \frac{C_{s,\kappa}}{1+|x|^{1+2s}},\end{equation}
for some~$C_{s,\kappa}>0$ depending only on~$s$ and~$\kappa$.

It is also useful to recall that, according to formula~(6) of~\cite{POLYA},
\begin{equation}\label{LYA}
\begin{split}&
\lim_{x\to\pm\infty} |x|^{1+2s}u(x,t)=
\lim_{x\to\pm\infty} |x|^{1+2s}
\int_{\R }
e^{-|2\pi \kappa\xi|^{2s}t}\cos(2\pi x \xi)
\,d\xi\\&\qquad=
2\lim_{x\to\pm\infty} |x|^{1+2s}
\int_0^{+\infty}
e^{-(2\pi\kappa\xi)^{2s}t}\cos(2\pi x \xi)
\,d\xi\\&\qquad=
\frac1{{\pi\kappa t^{\frac1{2s}}}}\lim_{x\to\pm\infty} |x|^{1+2s}
\int_0^{+\infty}
e^{-\vartheta^{2s}}\cos\left( \frac{ x \vartheta}{\kappa t^{\frac1{2s}}}\right)
\,{d\vartheta}\\&\qquad
=\frac{\kappa^{2s}\, t }{{\pi} }\lim_{y\to\pm\infty} |y|^{1+2s}
\int_0^{+\infty}
e^{-\vartheta^{2s}}\cos(y \vartheta )
\,{d\vartheta}\\&\qquad =\frac{\kappa^{2s}\,t\;\Gamma(1+2s)\,\sin(\pi s) }{{\pi} },
\end{split}
\end{equation}
where the substitutions~$\vartheta:=2\pi\kappa\xi t^{\frac1{2s}}$
and~$y:=\frac{ x }{\kappa t^{\frac1{2s}}}$
have been used.

\section{Description of the optimal strategies in different frameworks}\label{DIOFFE-R}

We introduce here the notion of value
functional related to the foraging success that we aim at optimizing
with respect to the parameter~$s$.

Given a distribution of targets~$p(x,t)$,
the {\em foraging success functional} will be taken as proportional to the encounters
between seekers and preys over time and therefore, given~$T>0$, it
takes the form
\begin{equation}\label{AV}
\iint_{\R\times(0,T)} p(x,t)\,u(x,t)\,dx\,dt
.\end{equation}
We will compare this quantity, which is advantageous for the forager, with several quantities
of interest which instead provide a penalization for the seeker's strategy. These terms will be {\em time}
(thus, we will consider the amount of targets met over the time span~$T$), a {\em renormalization of time}
that takes into account, in some sense, the trajectory performed at a discrete level by a corresponding L\'evy walker
(as presented in~\eqref{STRUC} below),
and the {\em average distance from the origin} (that is the distance of the forager ``from home'',
as discussed in~\eqref{AVEX} below).

To present the renormalization of time, we let~$s\in\left(\frac12,1\right)$
and we recall (see e.g. formula~(4.6) in~\cite{ABA}) that the mean excursion for each time step
of a discrete L\'evy walker is proportional to the spacial step by a factor of the form
\begin{equation*}
\frac{\displaystyle\sum_{j=1}^{+\infty} \frac{1}{j^{2s}}}{
\displaystyle\sum_{j=1}^{+\infty} \frac{1}{j^{1+2s}}}=\frac{\zeta(2s)}{\zeta(1+2s)}.
\end{equation*}
Though one cannot really consider this as the distance traveled by the L\'evy walker
in the unit of time (due to the nonlinear dependence between space and time variables
in long-jump random processes), it is suggestive to consider a possible renormalization of time of the form
\begin{equation}\label{STRUC} \overline\ell(s,T):=
\frac{T\,\zeta(2s)}{\zeta(1+2s)}.
\end{equation}

As for the distance between the forager and its burrow (located at the origin), we consider
the average displacement for~$s\in\left(\frac12,1\right)$ given by
\begin{equation}\label{AVEX} \ell(s,T):=\iint_{\R\times(0,T)} |x|\,u(x,t)\,dx\,dt.\end{equation}
We observe that this is a natural quantity to take into consideration as a penalization
for long excursions to account for the forager's need to return to home. Related
(but different) displacement functions were taken into account in equation~(1)
of~\cite{VIRTU}.
A variant of this approach (that will be accounted for in a forthcoming
work) consists in considering
pseudo mean displacements as in equation~(10)
of~\cite{RAP}, possibly also including
different normalization constants.

We also recall that the quantity in~\eqref{AVEX}
can be computed by using the Fourier Transform for generalized functions
(see Section~3.3 in Chapter~II of~\cite{GEL}
for the main results on this topic and Section~3.9 in Chapter~I
of~\cite{GEL}
for the setting of the notation related to generalized functions). 
Indeed, from\footnote{We stress that
the notation of~\cite{GEL} for the Fourier Transform chooses a different normalization than the one here,
by defining
$$ \widetilde{\mathcal{F}}f(\xi):=\int_{\R }f(x)e^{ix \xi}\,dx=
{\mathcal{F}}f\left(-\frac\xi{2\pi}\right),$$
see formula~(1) on page~153 in~\cite{GEL}.} equation~(2) on page~194 of~\cite{GEL} we know that
$$ {\mathcal{F}}(|x|)=-\frac1{2\pi^2\,|\xi|^2}$$
and therefore, by Plancherel Theorem, \eqref{WIDE}
and the substitution~$y:=(2\pi \kappa\xi)^{2s}t$,
\begin{eqnarray*}&&
\iint_{\R\times(0,T)} |x|\,u(x,t)\,dx\,dt=
\iint_{\R\times(0,T)} |x|\,\big(u(x,t)-\delta_0(x)\big)\,dx\,dt\\&&\qquad=-\frac1{2\pi^2}
\iint_{\R\times(0,T)} \frac{\widehat u(\xi,t)-1}{|\xi|^2}\,d\xi\,dt=
-\frac1{2\pi^2}
\iint_{\R\times(0,T)} \frac{
\exp\big(-|2\pi \kappa\xi|^{2s}t\big)
-1}{|\xi|^2}\,d\xi\,dt\\&&\qquad=
-\frac1{\pi^2}
\iint_{(0,+\infty)\times(0,T)} \frac{
\exp\big(-(2\pi \kappa\xi)^{2s}t\big)
-1}{\xi^2}\,d\xi\,dt\\&&\qquad=
-\frac{\kappa }{\pi s}
\iint_{(0,+\infty)\times(0,T)} t^{\frac1{2s}} \big(
e^{-y}
-1\big) y^{-\frac1{2s}-1}\,dy\,dt.
\end{eqnarray*}
In this way we obtain that
\begin{equation}\label{AVEX-2}\begin{split} \ell(s,T)\,&=\,
-\frac{2\kappa \,T^{\frac{1+2s}{2s}}}{\pi \,(1+2s)}
\int_0^{ +\infty } \big(
e^{-y}
-1\big) y^{-\frac1{2s}-1}\,dy\\&=
\frac{4\kappa s \,T^{\frac{1+2s}{2s}}}{\pi \,(1+2s)}
\int_0^{ +\infty } 
\left[ \frac{d}{dy}\left(
\frac{e^{-y}-1}{ y^{\frac1{2s}}}\right)+\frac{e^{-y}}{ y^{\frac1{2s}}}
\right]\,dy
\\&=
\frac{4\kappa s \,T^{\frac{1+2s}{2s}}}{\pi \,(1+2s)}
\int_0^{ +\infty } 
e^{-y}\, y^{\frac{2s-1}{2s}-1}\,dy\\&=
\frac{4\kappa s \,T^{\frac{1+2s}{2s}}}{\pi \,(1+2s)}
\Gamma\left(\frac{2s-1}{2s}\right).
\end{split}
\end{equation}

One of the main goals of this paper is to consider, as efficiency functional for the forager,
the ratio between~\eqref{AV} and either the time~$T$, or the quantity in~\eqref{STRUC},
or that in~\eqref{AVEX}. We stress that while~$T$ is obviously well defined\footnote{In any case, in what follows
we will restrict our analysis to the range~$s\in\left(\frac12,1\right)$.
This is due to the fact that the utility functions that we will define in~\eqref{UTIFU}, and thus in~\eqref{MK:03},
are properly defined only in this range (the functionals in~\eqref{LFT} could instead more generally be defined
in a larger range of~$s$).}
for all~$s\in(0,1)$, the
quantities in~\eqref{STRUC} and~\eqref{AVEX} are finite only when~$s\in\left(\frac12,1\right)$
(formally, they can be defined to be equal to~$+\infty$ when~$s\in\left(0,\frac12\right]$).
The reduction of the analysis of foraging strategies in the range~$s\in\left(\frac12,1\right)$
has been also performed elsewhere in the literature, see e.g. the discussion after
formula~(2)
in~\cite{Palyulin2931} or formula~(35) in~\cite{CHAOS}
(it is however interesting to pursue also different approaches to incorporate
conveniently modified situations in which the average jump distance is infinite,
but possibly incorporating waiting times between subsequent jumps,
see e.g.~\cite{ST1, ST2, SUG}
and pages~34--35 in~\cite{STRANGE}).

We will also distinguish two cases of interest according to the diffusion coefficient~$\kappa$ in~\eqref{BC:EQ}.
Namely, we will consider the standard case in which~$\kappa=1$
(this is a classical normalization choice, see e.g. formula~(1) in~\cite{Palyulin2931}), as well as the case in which~$\kappa$
depends on~$s$ via the relation
\begin{equation}\label{KAPS} \kappa=
\left( -\frac{ \cos(\pi s)\Gamma(-2s)}{\zeta(1+2s)}\right)^{\frac1{2s}}=:\kappa_s.\end{equation}
This form of the diffusion coefficient is the one emerging in the formal passage to the continuous
limit of a random L\'evy walker in the discrete lattice~$h\Z$ for time steps~$\tau=h^{2s}$, since, in this setting,
\begin{eqnarray*}
\frac{u(x,t+\tau)-u(x,t)}{\tau}&=&
\left( \sum_{k\in\Z\setminus\{0\}}\frac{1}{|k|^{1+2s}}\right)^{-1}\sum_{k\in\Z\setminus\{0\}}\frac{u(x+hk,t)-u(x,t)}{h^{2s} |k|^{1+2s}}
\\&=&
\left( 2\sum_{k=1}^{+\infty}\frac{1}{k^{1+2s}}\right)^{-1}\sum_{k=1}^{+\infty}\frac{u(x+hk,t)+u(x-hk,t)-2u(x,t)}{h^{2s} k^{1+2s}}\\
&\simeq&
\frac{1}{2\zeta(1+2s)}\int_0^{+\infty}\frac{u(x+y,t)+u(x-y,t)-2u(x,t)}{y^{1+2s}}\,dy,
\end{eqnarray*}
where we have approximated a Riemann sum with the corresponding integral. Thus setting~$v_{y,t}(x):=u(x+y,t)$, and noticing that
\begin{eqnarray*}
\widehat v_{y,t}(\xi)=
\int_{\R }u(x+y,t) e^{-2\pi ix \xi}\,dx=
e^{2\pi i y \xi}
\int_{\R }u(z,t) e^{-2\pi iz \xi}\,dz=e^{2\pi i y \xi}\,\widehat u(\xi,t),
\end{eqnarray*}
we see that, in the formal limit,
\begin{eqnarray*}
\partial_t u(x,t)
&=&
\frac{1}{2\zeta(1+2s)}
{\mathcal{F}}^{-1}\left(
\int_0^{+\infty}\frac{e^{2\pi i y \xi}+e^{-2\pi i y \xi}-2}{y^{1+2s}}\,dy\;
\widehat u(\xi,t)\right)\\
&=&
-\frac{1}{\zeta(1+2s)}
{\mathcal{F}}^{-1}\left(
\int_0^{+\infty}\frac{1-\cos(2\pi y\xi)}{y^{1+2s}}\,dy\;
\widehat u(\xi,t)\right)\\&=&
-\frac{(2\pi)^{2s}}{\zeta(1+2s)}\int_0^{+\infty}\frac{1-\cos z}{z^{1+2s}}\,dz\;
{\mathcal{F}}^{-1}\left(
|\xi|^{2s}\,
\widehat u(\xi,t)\right)\\&=&
\frac{ \cos(\pi s)\Gamma(-2s)}{\zeta(1+2s)}{\mathcal{F}}^{-1}\left(
|2\pi\xi|^{2s}\,
\widehat u(\xi,t)\right)\\&=&
\frac{ \cos(\pi s)\Gamma(-2s)}{\zeta(1+2s)}(-\Delta)^s u(x,t)
\end{eqnarray*}
see e.g. the appendix in~\cite{COCL} for the computation of the latter constant (which is negative),
and this justifies~\eqref{KAPS}. 

Moreover, using the functional equation~(40.5) in~\cite{ANT} (and, as customary,
adopting the notation that extends the Riemann Zeta Function by analytic continuation),
we can simplify the expression for $\kappa_s$ in~\eqref{KAPS} and get 
\begin{equation}\label{KAPPAS1.0}
\kappa_s=\frac{1}{2\pi}\left(-\frac{1}{2\zeta(-2s)}\right)^{\frac{1}{2s}}.
\end{equation}

\subsection{Single prey at the origin}\label{SINOR}

We now consider the case of a single target located at the origin.
In this case, the distribution of prey can be written as
$$ p_0(x)=\delta_0(x).$$
We observe that, by~\eqref{WIDE},
\begin{equation}\label{LB1}\begin{split}& \int_{\R}p_0(x)u(x,t)\,dx=u(0,t)={\mathcal{F}}^{-1}\big(
\widehat u(\cdot,t)\big)(0,t)
=\int_\R
\exp\big(-|2\pi \kappa\xi|^{2s}t\big)\,d\xi\\&\qquad\qquad=
2\int_{0}^{+\infty} \exp\left(- (2\pi\kappa \xi)^{2s}t\right)\,d\xi.\end{split}\end{equation}
Also, making use
of the change of variable~$\vartheta:=(2\pi \kappa\xi)^{2s}t$, we see that
\begin{equation}\label{LB} 
2\int_{0}^{+\infty} \exp\left(- (2\pi\kappa \xi)^{2s}t\right)\,d\xi=\frac{1}{2 \pi\kappa\, s\, t^{ \frac{1}{2s}}}
\int_{0}^{+\infty} \vartheta^{ \frac{1}{2s}-1} e^{-\vartheta}\,d\vartheta=\frac{1}{2\pi \kappa\,s \, t^{ \frac{1}{2s}}}\;
\Gamma\left( \frac{1}{2s}\right)
.\end{equation}
Thus, 
in the notation of~\eqref{AV}, using~\eqref{LB1} and~\eqref{LB},
the
foraging success functional for a single target located at the origin
takes the form, for~$s\in\left(\frac12,1\right)$,
\begin{equation}\label{and14}
\Phi_0(s;\kappa,T):=\iint_{\R\times(0,T)}p_0(x)u(x,t)\,dx\,dt
=\int_0^T
\frac{1}{2\pi \kappa\,s \, t^{ \frac{1}{2s}}}\;
\Gamma\left( \frac{1}{2s}\right)\,dt=
\frac{T^{\frac{2s-1}{2s}}}{\pi \kappa \,(2s-1)}
\Gamma\left( \frac{1}{2s}\right)
,\end{equation}
and takes value equal to~$+\infty$ when~$s\in\left(0,\frac12\right]$.
Hence, recalling~\eqref{STRUC}, \eqref{AVEX-2} and~\eqref{KAPPAS1.0},
we consider the utility functionals defined for~$s\in\left(\frac12,1\right)$
given by
\begin{equation}\label{UTIFU}
\begin{split}&{\mathcal{E}}_1(s;T):=\frac{\Phi_0(s;1,T)}{T}=
\frac{1}{\pi \,T^{\frac{1}{2s}}\,(2s-1)}
\Gamma\left( \frac{1}{2s}\right),\\
&{\mathcal{E}}_2(s;T):=\frac{\Phi_0(s;\kappa_s,T)}{T}=
\frac{2\,\big(-2\zeta(-2s)\big)^{\frac{1}{2s}}}{T^{\frac{1}{2s}}\,(2s-1)}
\Gamma\left( \frac{1}{2s}\right),\\
&{\mathcal{E}}_3(s;T):=\frac{\Phi_0(s;1,T)}{\overline\ell(s,T)}
=\frac{\zeta(1+2s)}{\pi \,T^{\frac{1}{2s}}\,(2s-1)\,\zeta(2s)}
\Gamma\left( \frac{1}{2s}\right),\\
&{\mathcal{E}}_4(s;T):=\frac{\Phi_0(s;\kappa_s,T)}{\overline\ell(s,T)}=
\frac{\zeta(1+2s)2\,\big(-2\zeta(-2s)\big)^{\frac{1}{2s}}
}{
 T^{\frac{1}{2s}}\,(2s-1)\,
\zeta(2s)}
\Gamma\left( \frac{1}{2s}\right),\\
&{\mathcal{E}}_5(s;T):=\frac{\Phi_0(s;1,T)}{\ell(s,T)}=
\frac{(1+2s)\;\Gamma\left( \frac{1}{2s}\right)}{4 s \,(2s-1)\,T^{\frac{1}{s}}\;
\Gamma\left(\frac{2s-1}{2s}\right)}
\\{\mbox{and }}\qquad
&{\mathcal{E}}_6(s;T):=\frac{\Phi_0(s;\kappa_s,T)}{ \ell(s,T)}=
\frac{\pi^2\big(-2\zeta(-2s)\big)^{\frac{1}{s}}\,(1+2s)\;\Gamma\left( \frac{1}{2s}\right)}{ s \,(2s-1)\,T^{\frac{1}{s}}\;
\Gamma\left(\frac{2s-1}{2s}\right)}.
\end{split}
\end{equation}

\begin{center}
\begin{figure}[!ht]
\includegraphics[width=0.40\textwidth]{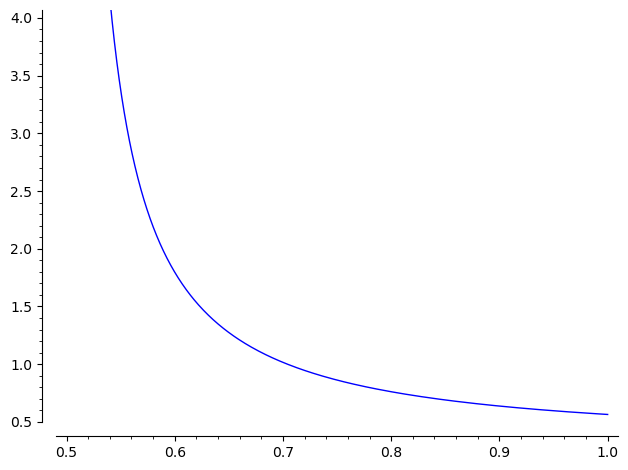}$\qquad$
\includegraphics[width=0.40\textwidth]{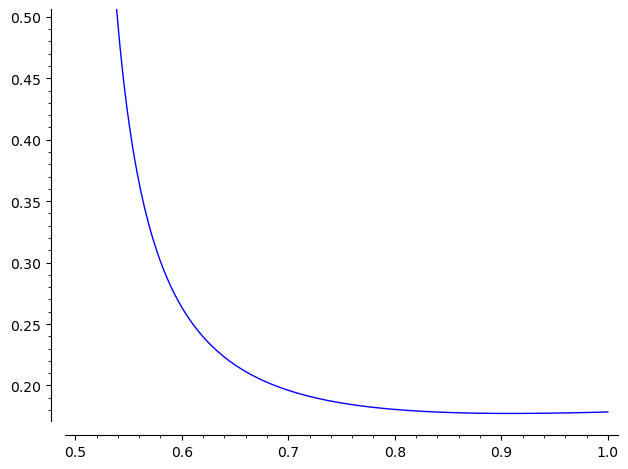}\\
\includegraphics[width=0.40\textwidth]{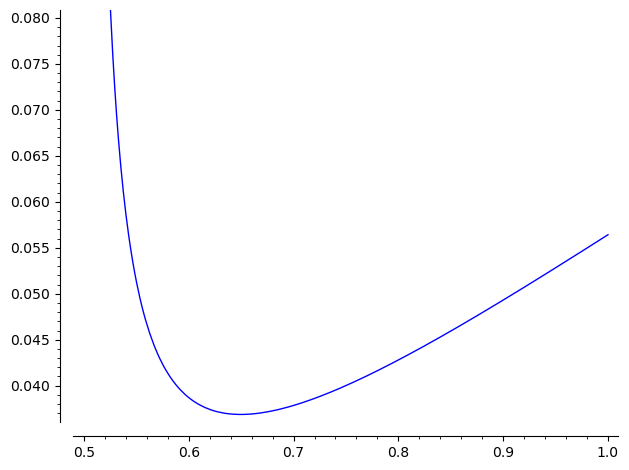}$\qquad$
\includegraphics[width=0.40\textwidth]{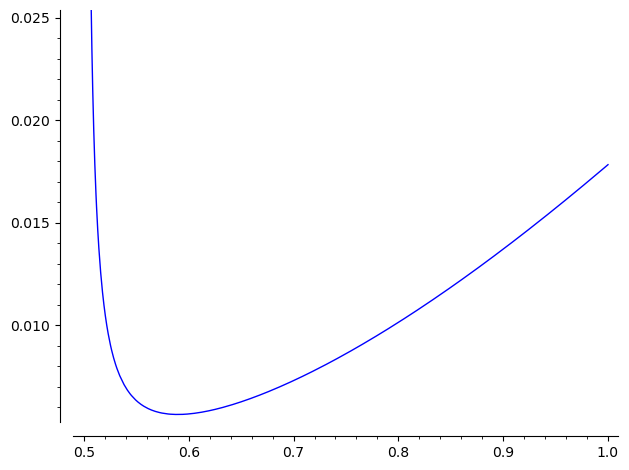}\\
\includegraphics[width=0.40\textwidth]{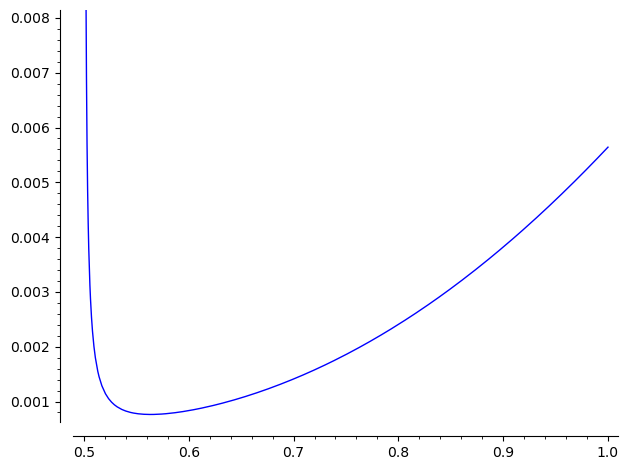}$\qquad$
\includegraphics[width=0.40\textwidth]{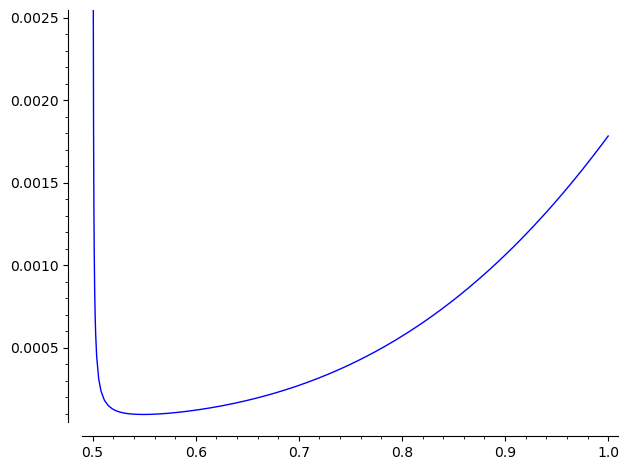}
\caption{\sl\footnotesize Plot of~$\left(\frac12,1\right)\ni s\mapsto
{\mathcal{E}}_1(s;T)$ for $T=10^j$, $j\in\{0,\dots,5\}$.}
        \label{63wDIK-MDyeufLO}
\end{figure}
\end{center}

We observe that
$$ \lim_{s\searrow1/2}{\mathcal{E}}_1(s;T)=+\infty,$$
therefore, for every~$T>0$.
\begin{equation}\label{E0d1}
{\mbox{the supremum of the utility functional }}
\left(\frac12,1\right)\mapsto {\mathcal{E}}_1(s;T){\mbox{
is uniquely attained at }}s=\frac12.
\end{equation}
We recall that the value~$s=\frac12$
occurs often in optimal foraging problems, as an ideal balance
between intensive search and longer (hence energetically more expensive) movements,
both
in terms of real world data (such as for atlantic cods, see e.g. Figure~1d in~\cite{SCA},
jackals, see e.g.~\cite{ATK}, 
wandering albatrosses,
see e.g. Figure~1 in~\cite{H2}, deers, see e.g. Figure~2(a)
in~\cite{DEE}, bees, see~\cite{BEE}, fruit flies, see~\cite{MOS}, 
and also Amazonian farmers searching for nuts, see Figure~3(b.09)
in~\cite{2018}, etc.)
and of theoretical optimization (see~\cite{VIS, BUL}). Interestingly,
it also occurs in patterns generated by human ecology (such as
distances between campsites, see
Figure~1 in~\cite{BRO}).
With respect to these data, the statement in~\eqref{E0d1}
can be seen as a confirmation of the most common paradigm
in the L\'evy foraging hypothesis. On the other hand,
the qualitative behaviour of~${\mathcal{E}}_1(s;T)$
changes dramatically for large intervals of time:
indeed, as hinted by Figure~\ref{63wDIK-MDyeufLO}
(that plots~${\mathcal{E}}_1(\cdot;T)$
for~$T\in\{1,10,10^2,10^3,10^4,10^5\}$), we have that
\begin{equation}\label{UNST}
\begin{split}&
{\mbox{for large~$T$, the utility functional }}
\left(\frac12,1\right)\mapsto {\mathcal{E}}_1(s;T)\\&{\mbox{has a unique minimum
at some point~$s_T$ such that}}\\
&\lim_{T\to+\infty} s_T=\frac12\qquad{\mbox{and}}\qquad
\lim_{T\to+\infty} 
{\mathcal{E}}_1(s_T;T)=0.
\end{split}
\end{equation}
This is an interesting phenomenon,
showing that the optimality at~$s=\frac12$ may become
``unstable'' and depends on the time span in which the phenomenon is
observed, allowing a sudden switch between
the optimal (but somehow unreliable) L\'evy foraging pattern
and the less ideal (but somehow more secure) classical Brownian motion
strategy.

It is suggestive to compare this phenomenon
to other occurrences in which
L\'evy flights
with~$s=\frac12$
should theoretically provide the optimal seeking strategy
but they coexist with another
possible notion of foraging optimization related to Brownian walks,
see e.g. the end of page~9 in~\cite{BUL}.

\begin{center}
\begin{figure}[h]
\includegraphics[width=0.55\textwidth]{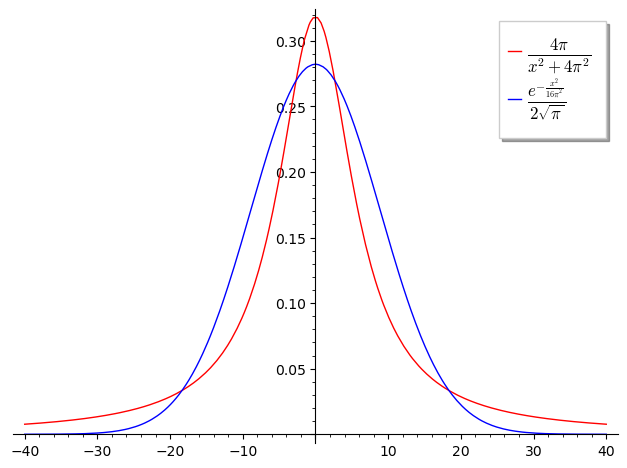}
\caption{\sl\footnotesize Plot of~$x\mapsto u(x,1)$ with~$\kappa=1$
corresponding to~$s=\frac12$ (in magenta)
and to~$s=1$ (in blue).}
        \label{6wyeufLOilb-ledue}
\end{figure}
\end{center}

A heuristic explanation for the statement in~\eqref{UNST}
can be given in terms of the behaviour
of the function~$u$ (with~$\kappa=1$)
at the origin and at infinity in dependence of the parameter~$s$.
Indeed, while the usual paradigm is to relate small values of~$s$
to long excursions of the traveller,
this general notion has sometimes to be revised according to the
specific mathematical model taken
into account in the 
diffusive strategy of the forager, since, on the one hand, solutions of
equation~\eqref{BC:EQ} corresponding to lower values of~$s$ do
present a fatter tail distribution,
but, on the other hand,
due to the loss of the regularizing effect of the diffusive operator for
small~$s$, they also present a
more prominent mass at the origin: see e.g. Figure~\ref{6wyeufLOilb-ledue}
in which one can compare solutions at time~$t=1$
corresponding to~$s=\frac12$ and~$s=1$.
Thus, the balance of these two apparently contrasting features
may provide advantageous foragers' strategies
for small values of~$s$ also in presence of proximate preys
(not due to the long range excursion induced
by the fat tail of the distribution, but rather due to the distribution peak
at the origin produced by the
less regularizing effect of a lower order operator).
With respect to this observation, in view of the scaling properties
of the equation (see~\eqref{SCALE}), the prominent role
of the peak at the origin occurs for small times, while it becomes
less significant for larger times. This somehow explains
why the L\'evy flights corresponding to~$s=\frac12$
are, in principle, more favorable than the classical Brownian motion,
but this effect may become less relevant and rather insecure
for very long time spans.
It is suggestive to investigate whether the interplay between
optimal but unstable strategies with suboptimal
but safer ones may play a role
in the appearance in nature of composite correlated random walks
and in the biological approximation of L\'evy walks as an innate
composite correlated random walks, see~\cite{MUSS0, MUSS1}.

It is also interesting to compare with
biological situations in which a predominance of classical random walks
coexists with patterns close to a theoretical optimum of~$s=\frac12$, see e.g. Figure~4 in~\cite{INTE}.

We stress that
the phenomenon described in~\eqref{UNST} relies on the ideal assumption that
the target is modelled as a ``material point'' (thus any arbitrarily small diffusion of the 
forager misses the resource) and is a byproduct of a memory-less
search strategy (see e.g.~\cite{MEMO} for a discussion of memory-enhanced foraging strategies).

The statement in~\eqref{UNST} can be checked analytically as detailed in Appendix~\ref{UNSTAPP}.

\begin{center}
\begin{figure}[h]
\includegraphics[width=0.40\textwidth]{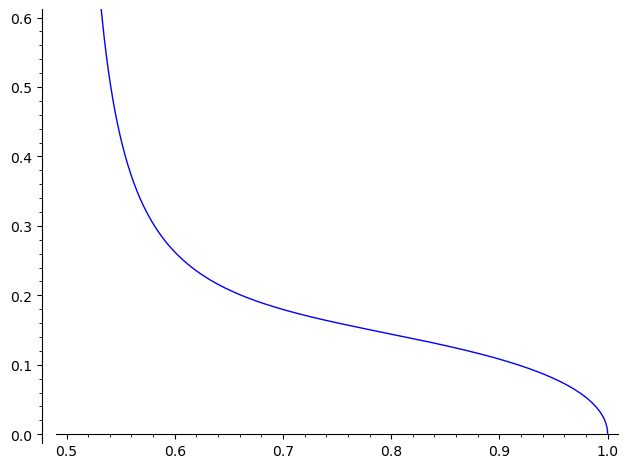}$\qquad$
\includegraphics[width=0.40\textwidth]{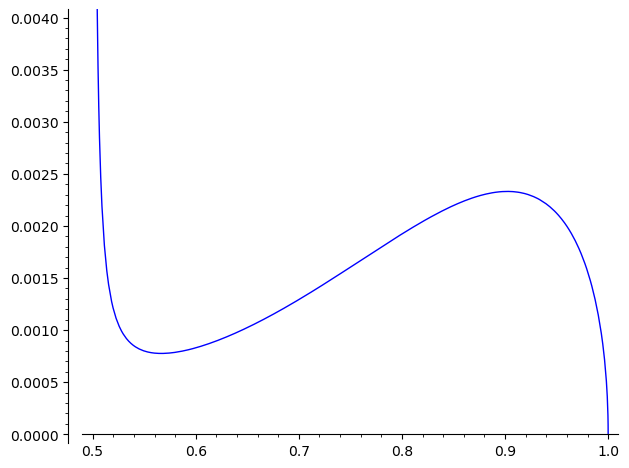}\\
\includegraphics[width=0.40\textwidth]{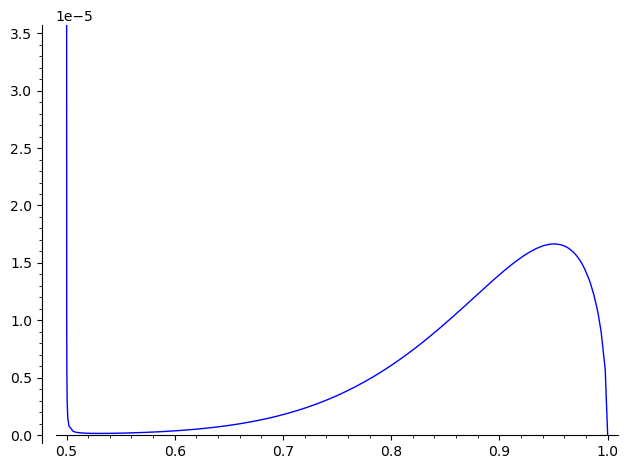}$\qquad$
\includegraphics[width=0.40\textwidth]{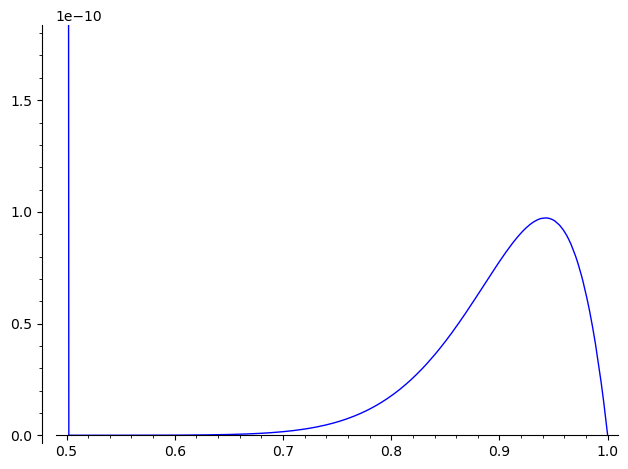}
\caption{\sl\footnotesize Plot of~$\left(\frac12,1\right)\ni s\mapsto
{\mathcal{E}}_2(s;T)$ for $T=10^{j}$, $j\in\{1,4,8,16\}$.}
        \label{63wDIK-MDyeufLO-2}
\end{figure}
\end{center}

We also have that
$$ \lim_{s\searrow1/2}{\mathcal{E}}_2(s;T)=+\infty$$
and, in view of the divergence of the Gamma Functions at negative integers,
$$ \lim_{s\nearrow1}{\mathcal{E}}_2(s;T)=0.$$
These equations show that, for every~$T>0$
\begin{equation}\label{E0d2}
\begin{split}&
{\mbox{the supremum of the utility functional }}
\left(\frac12,1\right)\mapsto {\mathcal{E}}_2(s;T){\mbox{
is uniquely attained at }}s=\frac12\\
&{\mbox{and the infimum is uniquely attained at }}s=1,\end{split}
\end{equation}
which in turn suggests a very strong advantage for the L\'evy
strategy compared to the Poisson one.
On the other hand, for long time spans, the pattern in~\eqref{E0d2}
shows a significant instability, as sketched in Figure~\ref{63wDIK-MDyeufLO-2}, which depicts
the map~${\mathcal{E}}_2(\cdot;T)$ for $T\in\{10,10^4,10^8,10^{16}\}$. As a result,
\begin{equation}\label{UNST-2}
\begin{split}&
{\mbox{for large~$T$, the utility functional }}
\left(\frac12,1\right)\mapsto {\mathcal{E}}_2(s;T)\\&{\mbox{has a local minimum
at some point~$s_T$ such that}}\\
&\lim_{T\to+\infty} s_T=\frac12\qquad{\mbox{and}}\qquad
\lim_{T\to+\infty} 
{\mathcal{E}}_2(s_T;T)=0,\\
&{\mbox{and a local maximum
at some point~$S_T$ such that } }
\lim_{T\to+\infty}S_T=1.\end{split}
\end{equation}
The statement in~\eqref{UNST-2} can be checked analytically as detailed in Appendix \ref{VF}.

\begin{center}
\begin{figure}[!ht]
\includegraphics[width=0.40\textwidth]{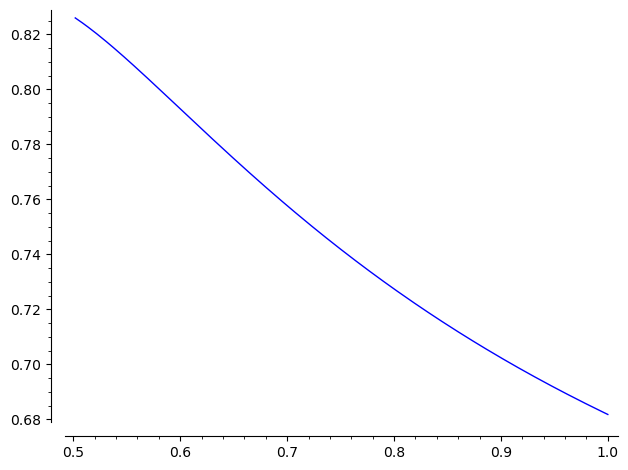}$\qquad$
\includegraphics[width=0.40\textwidth]{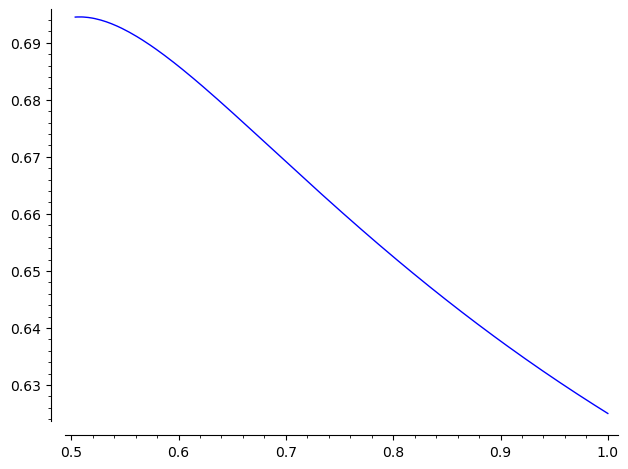}\\
\includegraphics[width=0.40\textwidth]{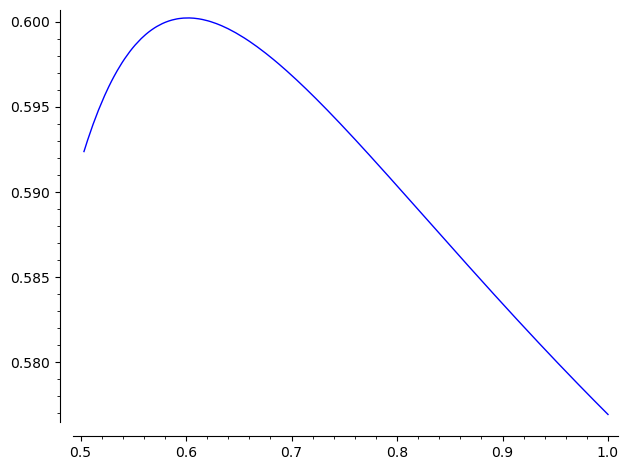}$\qquad$
\includegraphics[width=0.40\textwidth]{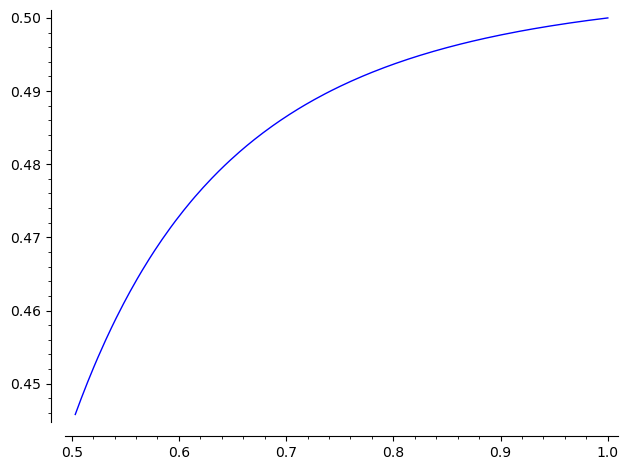}\\
\includegraphics[width=0.40\textwidth]{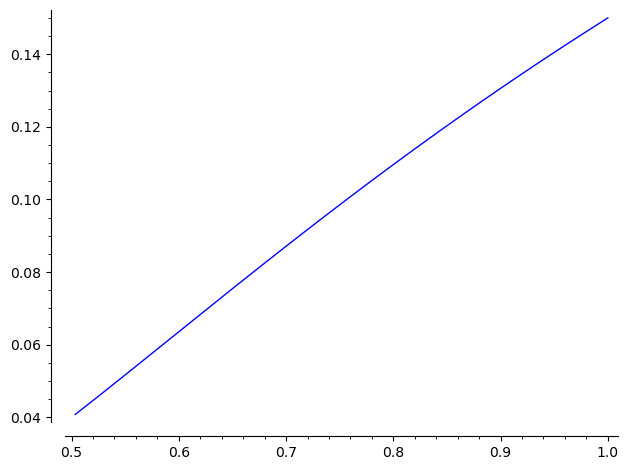}$\qquad$
\includegraphics[width=0.40\textwidth]{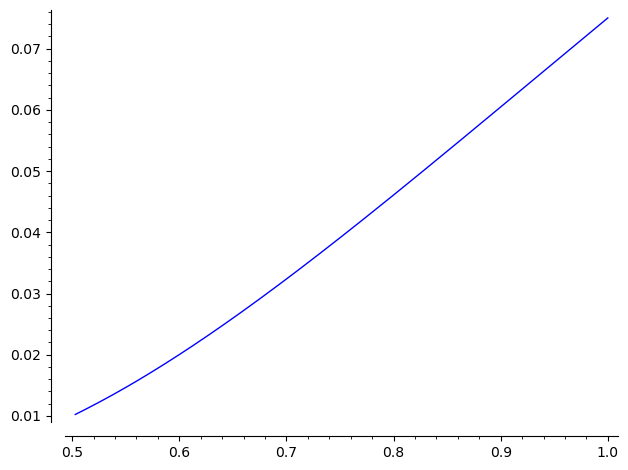}
\caption{\sl\footnotesize Plot of~$\left(\frac12,1\right)\ni s\mapsto
{\mathcal{E}}_5(s;T)$ for $T\in\{1.1, 1.2, 1.3, 1.5,5,10\}$.}
        \label{63wDIK-MDyeufLO-E5}
\end{figure}
\end{center}

\begin{center}
\begin{figure}[!ht]
\includegraphics[width=0.40\textwidth]{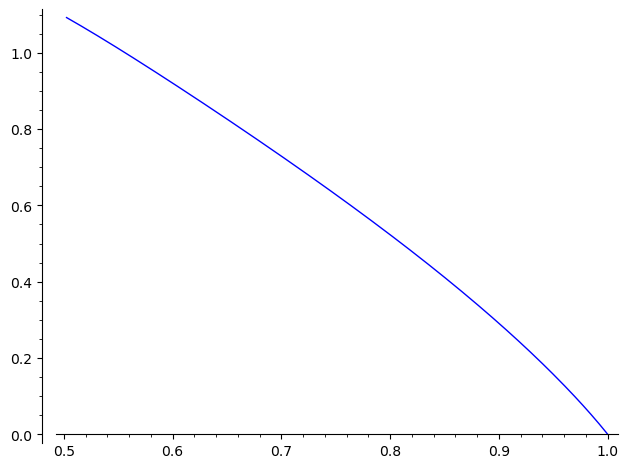}$\qquad$
\includegraphics[width=0.40\textwidth]{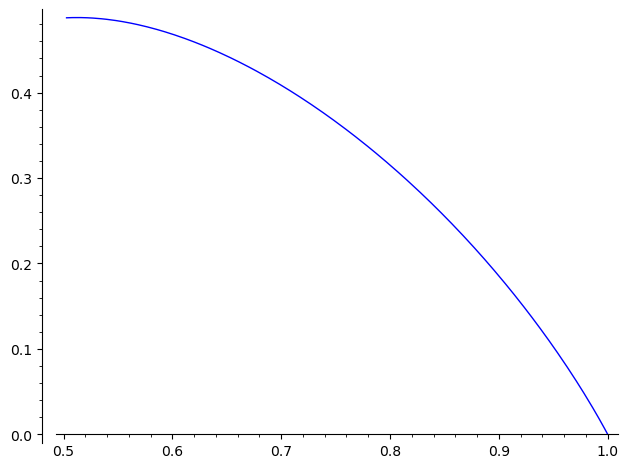}\\
\includegraphics[width=0.40\textwidth]{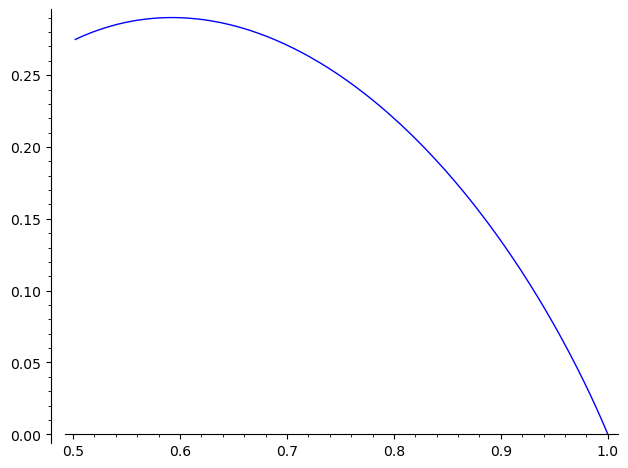}$\qquad$
\includegraphics[width=0.40\textwidth]{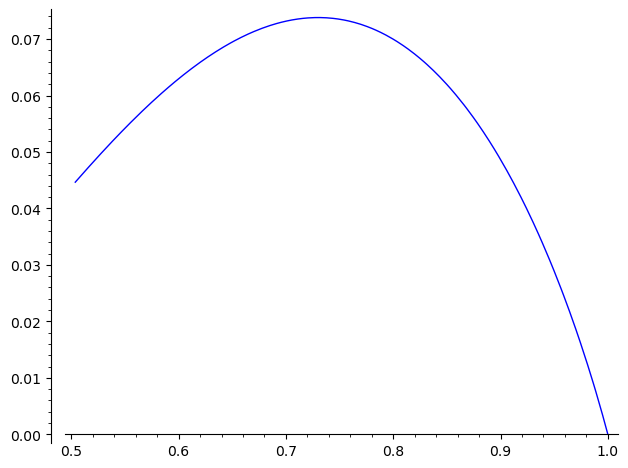}\\
\includegraphics[width=0.40\textwidth]{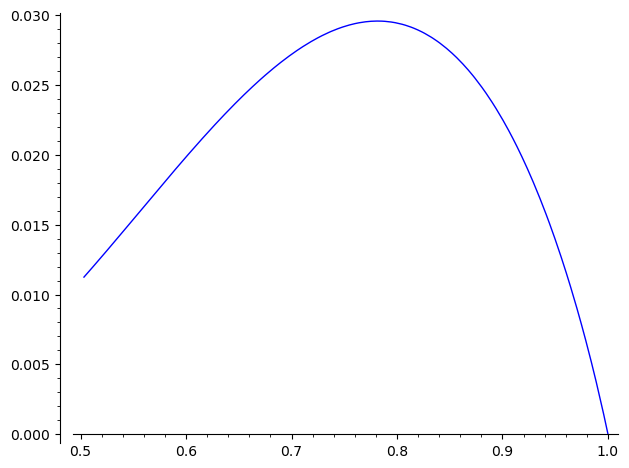}$\qquad$
\includegraphics[width=0.40\textwidth]{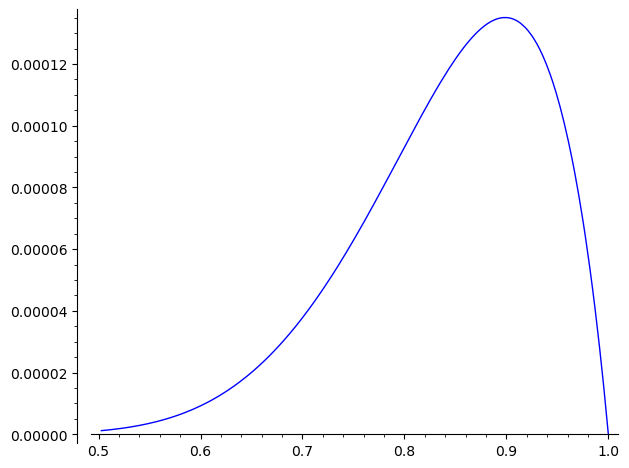}\\
\includegraphics[width=0.40\textwidth]{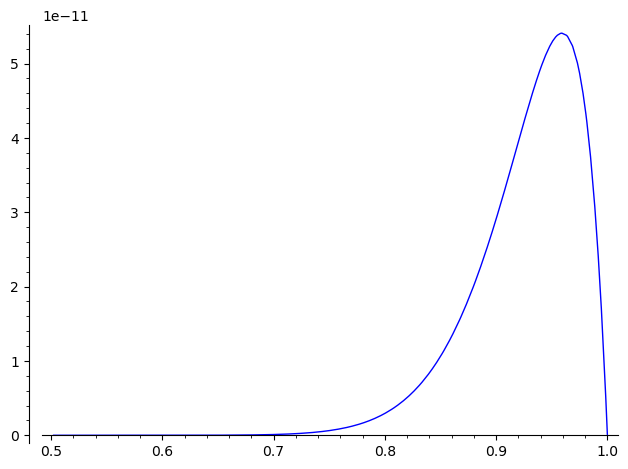}$\qquad$
\includegraphics[width=0.40\textwidth]{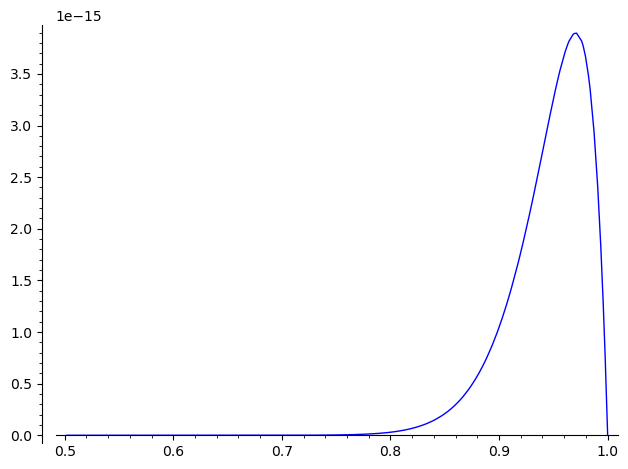}
\caption{\sl\footnotesize Plot of~$\left(\frac12,1\right)\ni s\mapsto
{\mathcal{E}}_6(s;T)$ for $T\in\{1,1.5,2,5,10,10^3,10^9,10^{13}\}$.}
        \label{63wDIK-MDyeufLO-E6}
\end{figure}
\end{center}

Interesting patterns 
having a theoretical (but unstable for large time) optimum at~$s=\frac12$
with a stabilizing option at~$s=1$
are exhibited by the utility functionals~${\mathcal{E}}_5$ and~${\mathcal{E}}_6$:
see Figure~\ref{63wDIK-MDyeufLO-E5}
for the sketch of~${\mathcal{E}}_5(\cdot;T)$
(notice the pattern change between~$T=1.1$ and~$T=1.5$
and the development of an interior maximum at~$T=1.3$)
and
Figure~\ref{63wDIK-MDyeufLO-E6}
for the sketch of~${\mathcal{E}}_6(\cdot;T)$.
We also stress that, for a given~$T$,
the values of~${\mathcal{E}}_5$ and~${\mathcal{E}}_6$
remain finite
(differently from the cases of~${\mathcal{E}}_1$ and~${\mathcal{E}}_2$).

\begin{center}
\begin{figure}[!ht]
\includegraphics[width=0.40\textwidth]{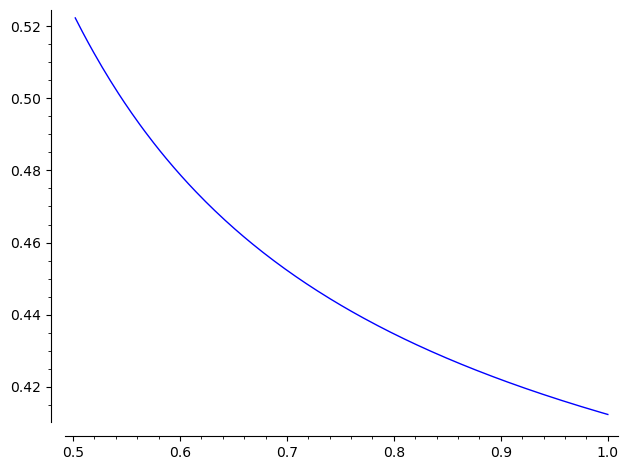}$\qquad$
\includegraphics[width=0.40\textwidth]{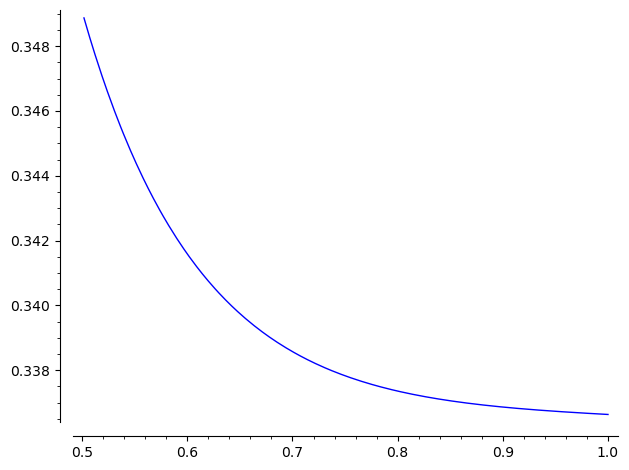}\\
\includegraphics[width=0.40\textwidth]{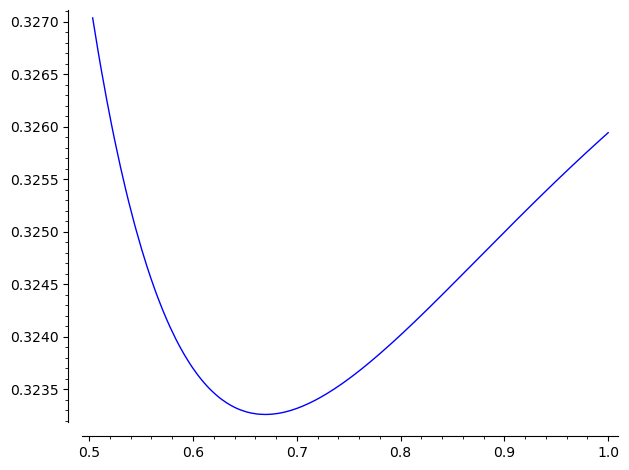}$\qquad$
\includegraphics[width=0.40\textwidth]{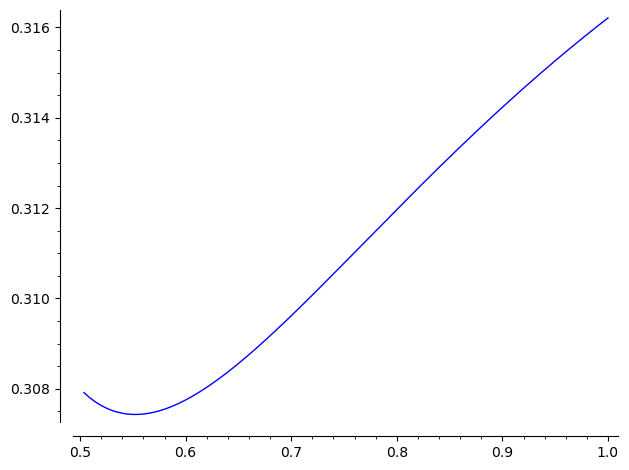}\\
\includegraphics[width=0.40\textwidth]{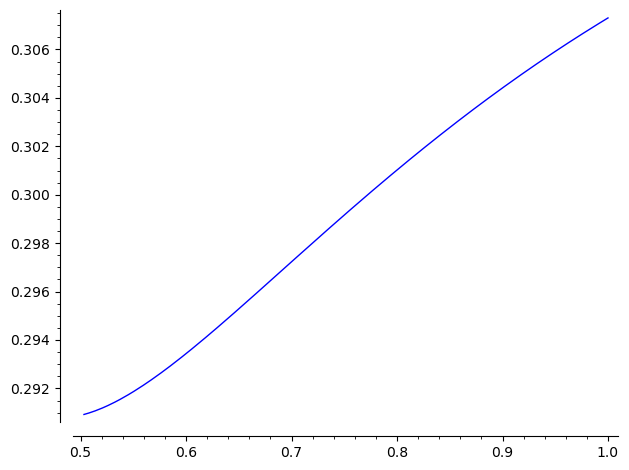}$\qquad$
\includegraphics[width=0.40\textwidth]{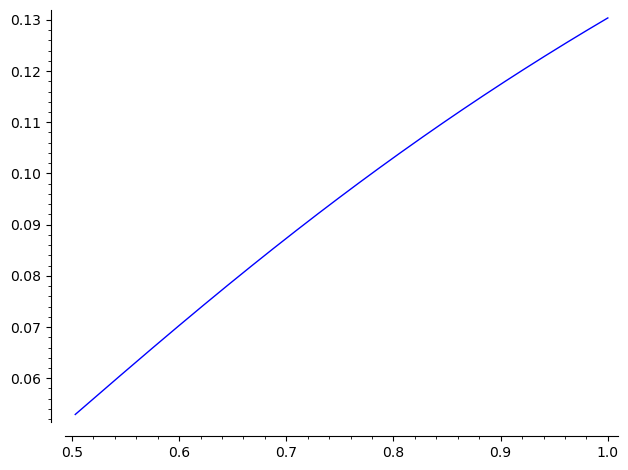}
\caption{\sl\footnotesize Plot of~$\left(\frac12,1\right)\ni s\mapsto
{\mathcal{E}}_3(s;T)$ for $T\in\{1,1.5,1.6,1.7,1.8,10\}$.}
        \label{63wDIK-MDyeufLO-Eil3}
\end{figure}
\end{center}

The cases of the utility functionals~${\mathcal{E}}_3$ and~${\mathcal{E}}_4$
are instead surprisingly different.
Indeed, Figure~\ref{63wDIK-MDyeufLO-Eil3}
hints that~${\mathcal{E}}_3$ is monotone decreasing with a supremum 
at~$s=\frac12$ when~$T\le 1.5$, but
then its monotonicity behavior changes when~$T\ge 1.6$
and develops a supremum 
at~$s=1$ when~$T\ge 1.7$. In this case,
even the theoretical optimality at~$s=\frac12$ is lost for large times
and additionally the switch between L\'evy and Poisson optimal strategy
occurs with a rather abrupt transition with respect to the parameter~$T$.

\begin{center}
\begin{figure}[!ht]
\includegraphics[width=0.40\textwidth]{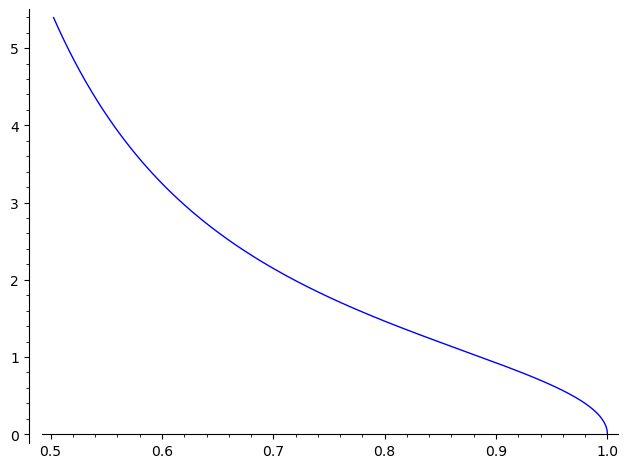}$\qquad$
\includegraphics[width=0.40\textwidth]{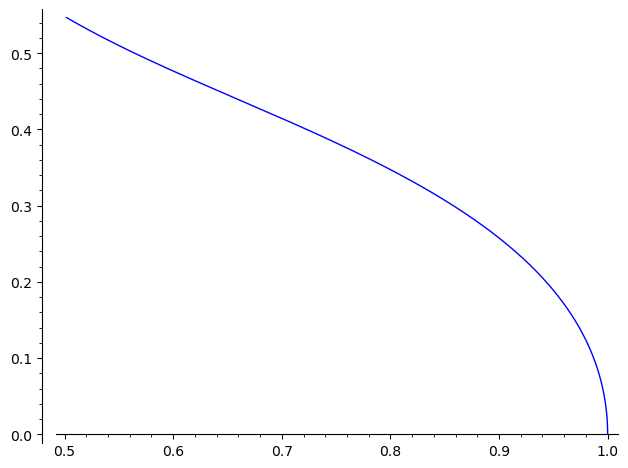}\\
\includegraphics[width=0.40\textwidth]{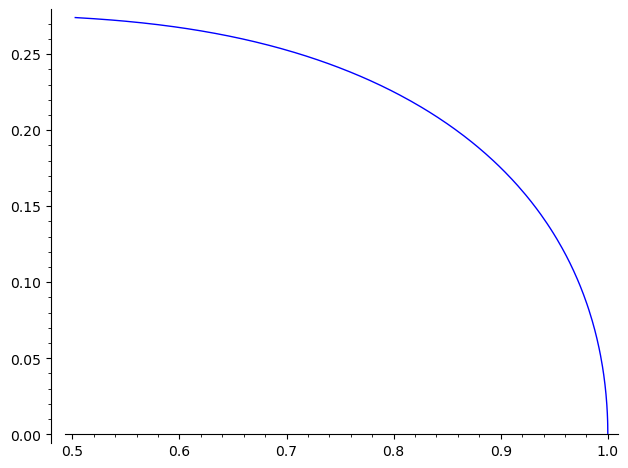}$\qquad$
\includegraphics[width=0.40\textwidth]{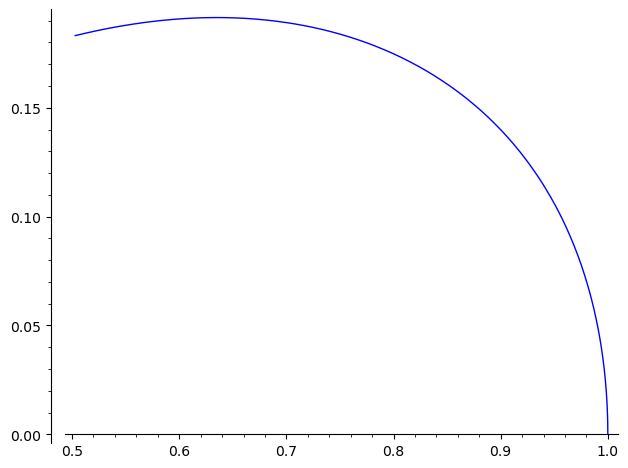}\\
\includegraphics[width=0.40\textwidth]{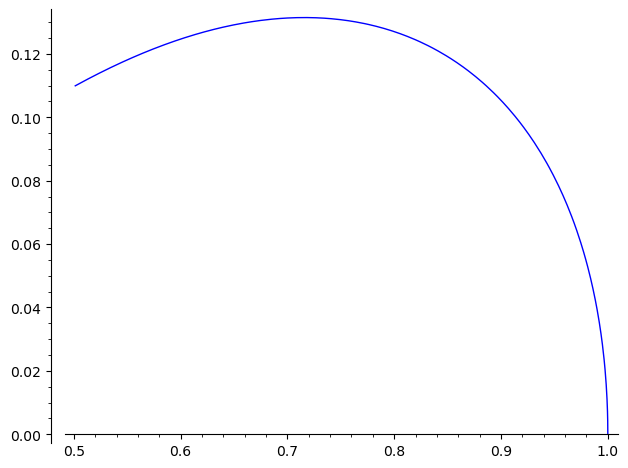}$\qquad$
\includegraphics[width=0.40\textwidth]{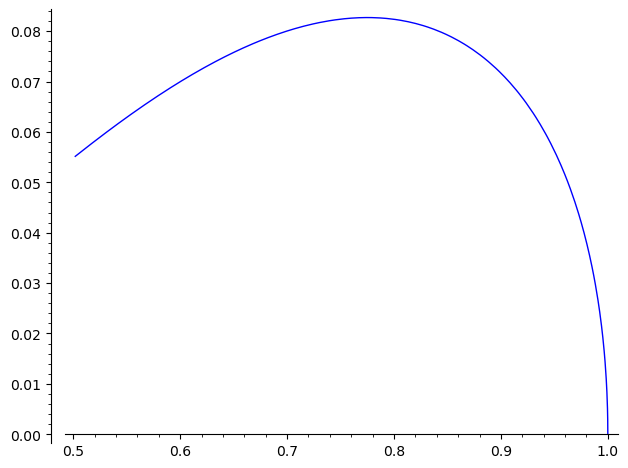}\\
\includegraphics[width=0.40\textwidth]{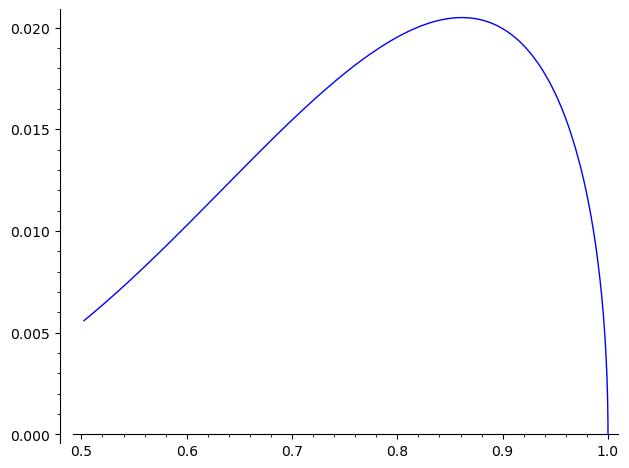}$\qquad$
\includegraphics[width=0.40\textwidth]{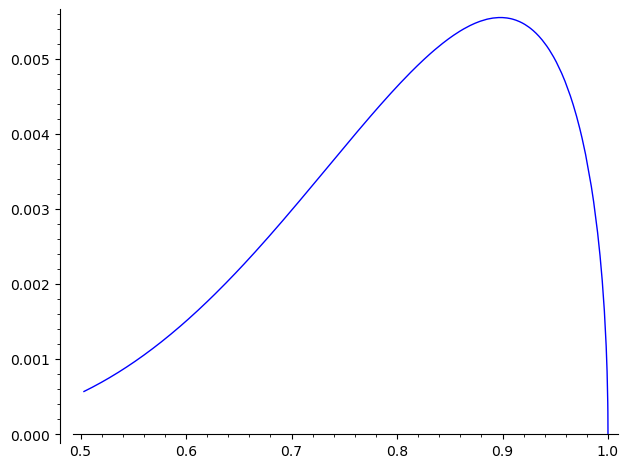}
\caption{\sl\footnotesize Plot of~$\left(\frac12,1\right)\ni s\mapsto
{\mathcal{E}}_4(s;T)$ for $T\in\{0.1, 1, 2, 3, 5, 10, 100, 1000\}$.}
        \label{63wDIK-MDyeufLO-Eil3OX}
\end{figure}
\end{center}

The functional~${\mathcal{E}}_4$ exhibits a different and interesting pattern,
as highlighted in Figure~\ref{63wDIK-MDyeufLO-Eil3OX}:
in this case the system shows
a sudden change of optimality occurring between~$T=2$ and~$T=3$:
it appears indeed that the L\'evy foraging for~$s=\frac12$ is optimal
when~$T\le 2$, but when~$T\ge3$ a new optimal strategy for a different fractional exponent~$s$ arises (with this new optimal exponent moving towards~$s=1$ as~$T$ becomes large).

As a matter of fact, we can detect analytically this bifurcation
phenomenon and find an explicit value for the critical~$T$ by
the following analytic argument.
We use the notation~$\epsilon:=2s-1$ according to which we have that
$$ (2s-1)\zeta(2s)=\epsilon\zeta(1+\epsilon)=
1 + \gamma\epsilon+o(\epsilon)$$
and therefore
\begin{eqnarray*}
{\mathcal{E}}_4(s;T)&=&
\frac{2\,\zeta(2+\epsilon)\;(-2\zeta(-\epsilon-1))^{\frac{1}{1+\epsilon}}}{T^{\frac{1}{1+\epsilon}}\;(1+\epsilon\gamma+o(\epsilon))}\,\Gamma\bigg(\frac{1}{1+\epsilon}\bigg)\\
&=&\frac{2\,(\zeta(2)+\zeta'(2)\epsilon+o(\epsilon))\,(-2\zeta(-1)+2\zeta'(-1)\epsilon+o(\epsilon))^{\frac{1}{1+\epsilon}}}{(T-T\ln T\epsilon+o(\epsilon))(1+\gamma\epsilon+o(\epsilon))}\,(1+\gamma\epsilon+o(\epsilon))\\
&=& 2\,\bigg(\frac{\pi^2}{6}+\zeta'(2)\epsilon+o(\epsilon)\bigg)\bigg(\frac{1}{6}+\frac{1}{6}\ln 6\,\epsilon+2\zeta'(-1)\epsilon+o(\epsilon)\bigg)\bigg(\frac{1}{T}+\frac{\ln T}{T}\epsilon+o(\epsilon)\bigg)+o(\epsilon)\\
&=&\frac{\pi^2}{18\,T}+2\epsilon\bigg(\frac{\pi^2}{36\,T}\ln T+\frac{\pi^2}{36\,T}\ln6+\frac{\pi^2}{3\,T}\zeta'(-1)+\frac{\zeta'(2)}{6\,T}\bigg)+o(\epsilon)
\end{eqnarray*}

This yields that
\begin{eqnarray*}
\frac{9 T}{\pi^2}\,\frac{d{\mathcal{E}}_4}{ds}\left(\frac12;T\right)&=&
 \ln T + \ln6+12 \zeta'(-1) + \frac{6}{\pi^2}\zeta'(2)
\end{eqnarray*}
and this quantity is positive (respectively, negative) when~$T> T_\star$
(respectively, when~$T< T_\star$), where
$$ T_\star:=
\exp\left(-\ln6-12\zeta'(-1)-\frac{6}{\pi^2}\zeta'(2)
\right)=2.145248182...
$$
and, as a result, when~$T>T_\star$ the inverse square law~$s=\frac12$
cannot maximize~${\mathcal{E}}_4$ 
and values of~$s$ even slightly larger than~$\frac12$ provide greater
values of such an efficiency functional.

\subsection{Preys on a sparse lattice}

Now we take into consideration a set of targets displayed in the lattice~$\lambda\Z$,
for some~$\lambda>0$, and we consider the asymptotics related to large values of~$\lambda$,
corresponding to the case of sparse preys. To this end, we consider the target distribution
\begin{equation*} p(x):=\sum_{k\in\Z } \delta_{\lambda k}(x).\end{equation*}
We observe that
\begin{equation}\label{pp} \int_{\R } p(x) u(x,t)\,dx =
\sum_{k\in\Z } u(\lambda k,t)=\sum_{k\in\Z } u_\lambda(k,t)
\end{equation}
where
$$ u_\lambda(x,t):=u(\lambda x,t).$$
By~\eqref{SCALE} and~\eqref{GIU},
\begin{equation}\label{PO}
\begin{split}&
(1+|x|^{1+2s}) |u_\lambda(x,t)|=
(1+|x|^{1+2s}) |u(\lambda x,t)|=
\frac{1+|x|^{1+2s}}{t^{ \frac{1}{2s}}}
u\left( \frac{\lambda x}{t^{\frac1{2s}}},1\right)\\&\qquad\le
\frac{C\,(1+|x|^{1+2s})}{t^{ \frac{1}{2s}}\left( 1+\left|\frac{\lambda x}{t^{\frac1{2s}}}\right|^{1+2s} \right)}\le
C_{s,\lambda,t}
,\end{split}
\end{equation} 
for some~$C_{s,\lambda,t}\in(0,+\infty)$.

Moreover, by~\eqref{WIDE},
\begin{eqnarray*}&& \widehat u_\lambda(\xi,t)=\int_{\R }u(\lambda x,t)e^{-2\pi ix \xi}\,dx=
\frac1{\lambda }
\int_{\R }u(y,t)e^{-2\pi i\lambda^{-1}y \xi}\,dy\\&&\quad\qquad\qquad=\frac1{\lambda }\widehat u\left(\frac\xi\lambda,t\right)
=\frac1{\lambda } \exp\left(-\frac{|2\pi \kappa\xi|^{2s}t}{\lambda^{2s}}\right)
\end{eqnarray*}
and accordingly
\begin{equation*}
\begin{split}&
(1+|\xi|^{1+2s}) |\widehat u_\lambda(\xi,t)|=
\frac{1+|\xi|^{1+2s}}{\lambda } \exp\left(-\frac{|2\pi \kappa\xi|^{2s}t}{\lambda^{2s}}\right)
\le
\widetilde C_{s,\lambda,t}
,\end{split}
\end{equation*} for some~$\widetilde C_{s,\lambda,t}\in(0,+\infty)$.

In view of this estimate and~\eqref{PO}, for a given~$t>0$
we can use the Poisson Summation Formula on~$u_\lambda$
(see e.g. formula~(4.4.2) and Theorem~4.4.2 in~\cite{PIN}) and, in light of~\eqref{pp}, conclude that
$$ \int_{\R } p(x)u(x,t)\,dx=\sum_{k\in\Z }u_\lambda(k)=
\sum_{k\in\Z }\widehat u_\lambda(k)=\frac1{\lambda } 
\sum_{k\in\Z } \exp\left(-\frac{|2\pi k|^{2s}t}{\lambda^{2s}}\right).
$$
For large~$\lambda$, we can consider the latter term as a Riemann sum, therefore,
using polar coordinates,
$$ \int_{\R } p(x)u(x,t)\,dx\simeq
\int_{\R } \exp\left(- |2\pi y|^{2s}t\right)\,dy=
2\int_{0}^{+\infty} \exp\left(- (2\pi \rho)^{2s}t\right)\,d\rho.
$$
Thus, recalling~\eqref{AV}, \eqref{LB} and~\eqref{and14},
we can consider, for large~$\lambda$, the foraging success functional
$$ \Phi(s):=\iint_{\R\times(0,T)}p(x)u(x,t)\,dx\,dt\simeq
\frac{1}{2\pi s }\;
\Gamma\left(\frac{1}{2s}\right)\int_0^T\frac{dt}{
t^{\frac{1}{2s}}}=\frac{T^{\frac{2s-1}{2s}} }{\pi(2s-1)}\;\Gamma\left(\frac{1}{2s}\right)=
\Phi_0(s;1,T).
$$
The case of a sparse distribution of targets is therefore reduced
to that of a single prey at the origin and, since the optimizers discussed
in Section~\ref{SINOR}
were isolated and nondegenerate,
the analysis provided in Section~\ref{SINOR} for a single prey gives asymptotic information
to the case of sparsely distributed targets
when~$\lambda$ is sufficiently large.

\subsection{Remote single prey}\label{REMOTE}

Now we consider the case of a single target located far away from the initial position of the seeker.
For this, given~$L>0$, let
$$ p_L(x):=\delta_L(x).$$
For~$T>0$, in view of~\eqref{AV}, we consider the foraging success functional
\begin{equation*}
\Phi_{L,T}(s;\kappa):=\iint_{\R\times(0,T)}p_L(x)u(x,t)\,dx\,dt=
\int_{0}^T u(L,t)\,dt
\end{equation*}
By~\eqref{LYA}, for large~$L$,
\begin{equation*}
L^{1+2s} \,\Phi_{L,T}(s;\kappa)\simeq
\int_{0}^T \frac{\kappa^{2s}\,t\;\Gamma(1+2s)\,\sin(\pi s) }{{\pi} }\,dt=
\frac{\kappa^{2s}\,T^2 \;\Gamma(1+2s)\,\sin(\pi s) }{{2}\pi }.
\end{equation*}
Thus, 
in the lines of~\eqref{UTIFU}, we set
\begin{equation*}
\widetilde\Phi_{L,T}(s;\kappa):=
\frac{\kappa^{2s}\,T^2 \;\Gamma(1+2s)\,\sin(\pi s) }{{2}\pi L^{1+2s}},
\end{equation*}
we notice that
\begin{equation}\label{APPEROS}
\Phi_{L,T}(s;\kappa)\simeq\widetilde\Phi_{L,T}(s;\kappa)
\end{equation}
for large~$L$, and we
discuss the optimization of 
the utility functionals\footnote{The final expressions for $\mathcal{G}_2$ and $\mathcal{G}_4$ are due to Euler's reflection formula
$$
\Gamma(1-z)\Gamma(z)=\frac{\pi}{\sin(\pi z)}
\qquad {\mbox{ for all }} z \in \mathbb{R}\setminus \mathbb{Z}.$$}
\begin{equation}\label{LFT}
\begin{split}&{\mathcal{G}}_1(s;L,T):=\frac{\widetilde\Phi_{L,T}(s;1)}{T}=
\frac{T \;\Gamma(1+2s)\,\sin(\pi s) }{{2}\pi \,L^{1+2s}},\\
&{\mathcal{G}}_2(s;L,T):=\frac{\widetilde\Phi_{L,T}(s;\kappa_s)}{T}=
-\frac{T \Gamma(-2s)\,\Gamma(1+2s)\,\sin(2\pi s) }{{4}\pi \,L^{1+2s}
\,\zeta(1+2s)}=
\frac{T}{{4}L^{1+2s}\, \zeta(1+2 s)}
 ,\\
&{\mathcal{G}}_3(s;L,T):=\frac{\widetilde\Phi_{L,T}(s;1)}{\overline\ell(s,T)}
=
\frac{T \;\zeta(1+2s)\,\Gamma(1+2s)\,\sin(\pi s) }{{2}\pi L^{1+2s}\,\zeta(2s)},\\
&{\mathcal{G}}_4(s;L,T):=\frac{\widetilde\Phi_{L,T}(s;\kappa_s)}{
\overline\ell(s,T)}=
-
\frac{T \,\Gamma(-2s)\,\Gamma(1+2s)\,\sin(2\pi s) }{{4}\pi L^{1+2s}\,\zeta(2s)}=
\frac{T }{{4} L^{1+2s}\zeta(2 s)},\\
&{\mathcal{G}}_5(s;L,T):=\frac{\widetilde\Phi_{L,T}(s;1)}{\ell(s,T)}=
\frac{T^{\frac{2s-1}{2s}} \;\Gamma(2+2s)\,\sin(\pi s) }{{8}L^{1+2s}\,\Gamma\left(\frac{2s-1}{2s}\right)}
\\{\mbox{and }}\qquad
&{\mathcal{G}}_6(s;L,T):=\frac{\widetilde\Phi_{L,T}(s;\kappa_s)}{ \ell(s,T)}=
\left(\frac{-1}{2^{1+2s}\pi^{2s}\zeta(-2s)}\right)^{\frac{2s-1}{2s}}
\frac{T^{\frac{2s-1}{2s}} \;\Gamma(2+2s)\,\sin(\pi s) }{{8}L^{1+2s}\,
\Gamma\left(\frac{2s-1}{2s}\right)}.
\end{split}\end{equation}
We point out that the final time~$T$ does not play
any role in the optimization in~$s$ of the value functionals~${\mathcal{G}}_1$,
${\mathcal{G}}_2$, ${\mathcal{G}}_3$ and~${\mathcal{G}}_4$
in~\eqref{LFT}.
We also stress that the biological meaning of
the efficiency functionals in~\eqref{LFT}
only occurs for large values of~$L>0$,
due to the asymptotics in~\eqref{APPEROS},
nevertheless it is interesting to
study those functionals for all values of~$L$
also to detect bifurcation phenomena with respect to this
parameter that depend only on the final analytic formulation
and not on their initial construction.

We also observe that, differently from the functionals in~\eqref{UTIFU}, the ones defined in~\eqref{LFT}
can be continuously extended in $\left(0,\frac{1}{2}\right]$. Nevertheless, among them only $\mathcal{G}_1$ and $\mathcal{G}_2$ admit non negative values for $s\in \left(0,\frac{1}{2}\right]$.
Since the functionals in~\eqref{LFT} lose their physical meaning for negative values,
when studying the optimal search strategy we will take into account the fractional parameter~$s$ in the whole
interval~$(0,1)$ only when studying $\mathcal{G}_1$ and $\mathcal{G}_2$.

\begin{center}
\begin{figure}[!ht]
\includegraphics[width=0.40\textwidth]{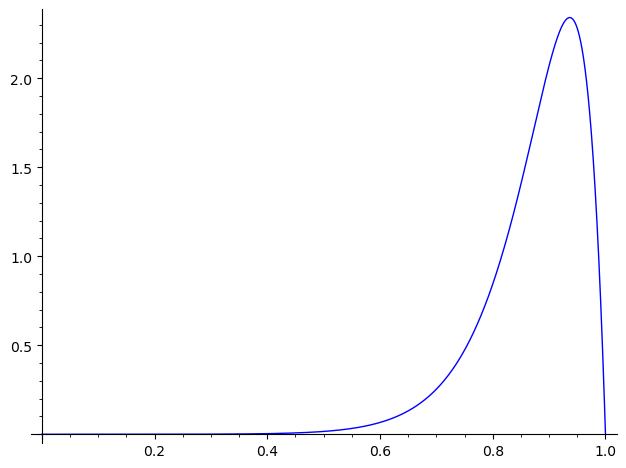}$\qquad$
\includegraphics[width=0.40\textwidth]{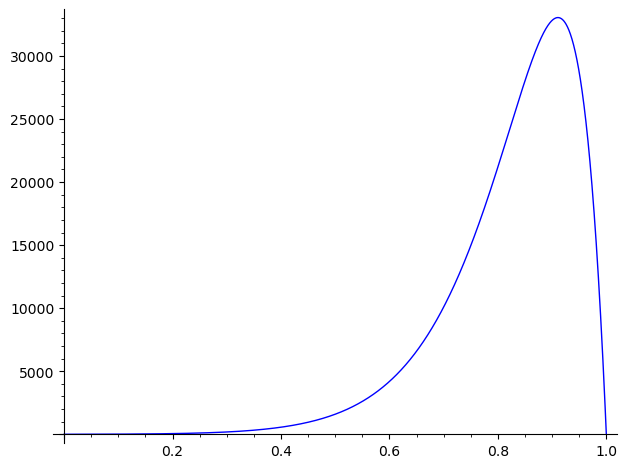}\\
\includegraphics[width=0.40\textwidth]{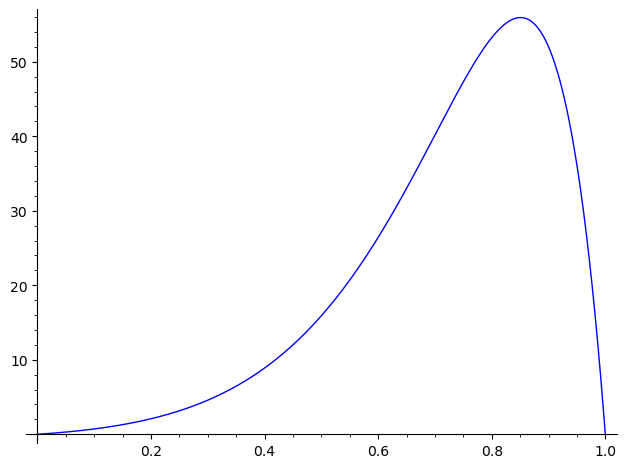}$\qquad$
\includegraphics[width=0.40\textwidth]{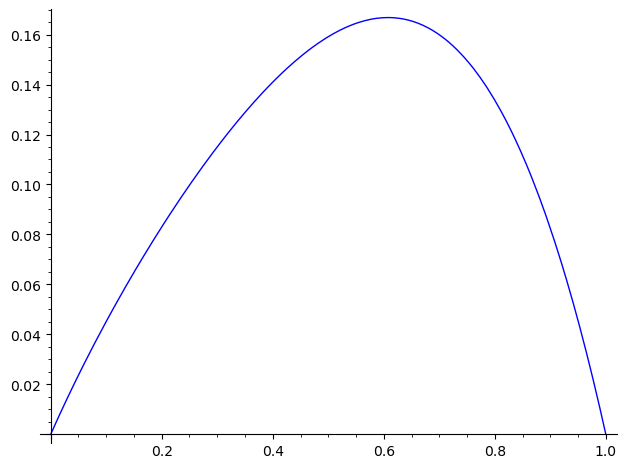}\\
\includegraphics[width=0.40\textwidth]{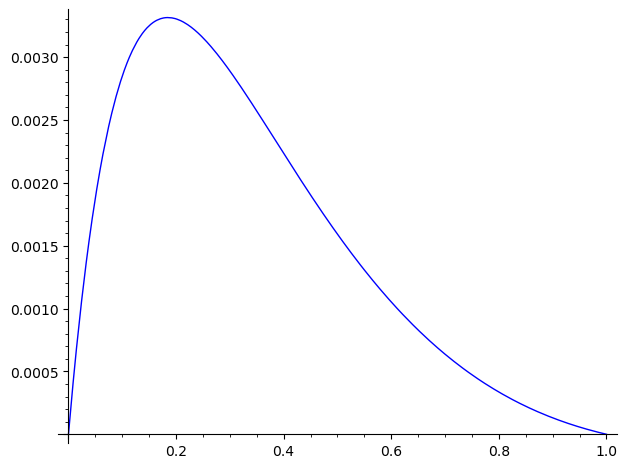}$\qquad$
\includegraphics[width=0.40\textwidth]{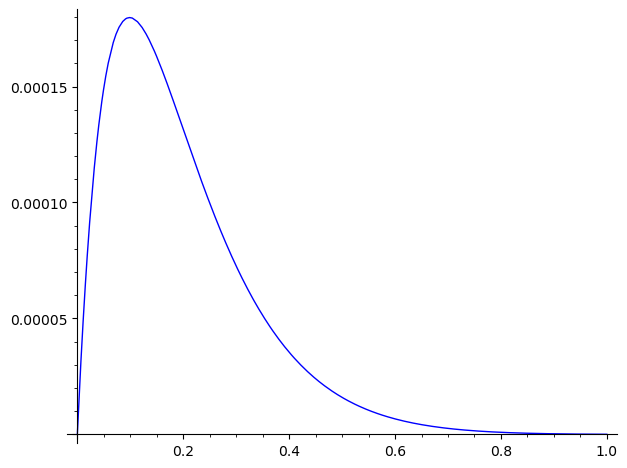}\\
\caption{\sl\footnotesize Plot of~$\left(0,1\right)\ni s\mapsto
{\mathcal{G}}_1(s;L,T)$ (say, when~$T=1$)
for~$L=10^j$, with~$j\in\{-3,\dots,2 \}$.}
        \label{63wDIK-MD89907-04yeufLO}
\end{figure}
\end{center}

A plot of~${\mathcal{G}}_1$ is given in Figure~\ref{63wDIK-MD89907-04yeufLO},
where one can appreciate that 
for small values of~$L$ the optimal exponent
is close (but not equal) to~$1$
(consistently with the idea that if the prey is
close to the starting point of the predator the best seeking
strategy is close to that of local type),
while for large values of~$L$ the maximum
is provided by a value of $s$ close to zero. This is in agreement with the  idea that if the
prey is far away from the forager's starting position, then the optimal search strategy has a strong nonlocal component. More precisely, we have that 
\begin{equation}\label{SLL-1}
\begin{split}&
{\mbox{if $L$ is large enough, the optimal foraging strategy
for~${\mathcal{G}}_1(s;L,T)$}}\\
&{\mbox{is uniquely attained at some~$s_L$ such that }}\lim_{L\to+\infty}
s_L=0.\end{split}\end{equation}
This can be proved analytically, see Appendix~\ref{ALLc-1}. From the plots provided in
Figure~\ref{63wDIK-MD89907-04yeufLO}, we see that $s=0$ is a global minimum for $\mathcal{G}_1$.
Nevertheless, as stated in \eqref{SLL-1}, for $L$ large enough $\mathcal{G}_1$ admits a unique
maximum~$s_L\in (0,1)$ such that~$s_L\searrow 0$ for~$L\to +\infty$. Thus, the global maximum~$s_L$
approaches the global minimum~$s=0$ as~$L$ gets larger and larger.
We can summarize this phenomenon, by saying that the optimal search strategy for $\mathcal{G}_1$ is unstable. Namely, a slight deviation from the optimal value $s_L$ can lead to very small values for $\mathcal{G}_1$. 

The functionals that we take into account to model the environmental scenario of a remote single prey are obtained using the approximation provided in~\eqref{APPEROS}. For this reason, we cannot rule out a priori that the aforementioned instability result is a consequence of this approximation and instead
it does occur if we consider the original functional 
\begin{equation}\label{akuybgefc6-1}
\frac{\Phi_{L,T}(s;1)}{T}.
\end{equation} 
In a forthcoming paper, see~\cite{FI}, we will address this problem and we will show that this instability
result does hold true also for the non approximated efficiency functional in~\eqref{akuybgefc6-1},
provided that $L$ is large enough. We will also show that this unstable behaviour is true in a multidimensional framework. 

\begin{center}
\begin{figure}[!ht]
\includegraphics[width=0.40\textwidth]{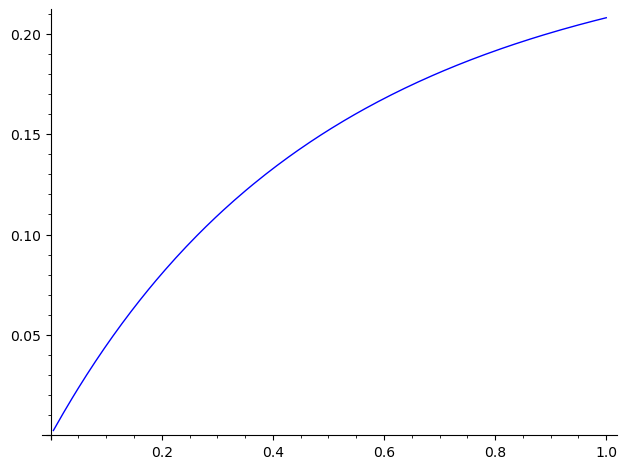}$\qquad$
\includegraphics[width=0.40\textwidth]{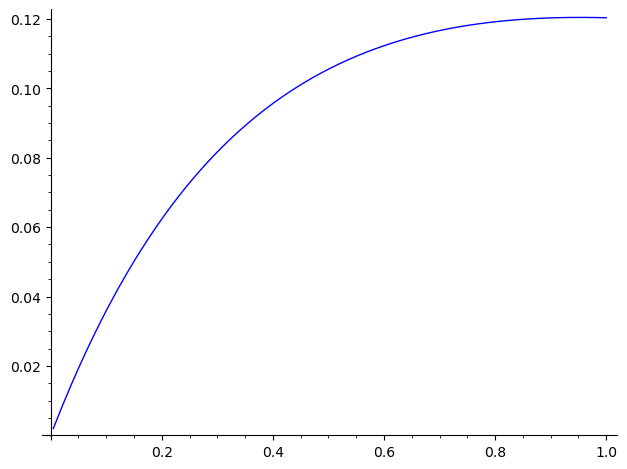}\\
\includegraphics[width=0.40\textwidth]{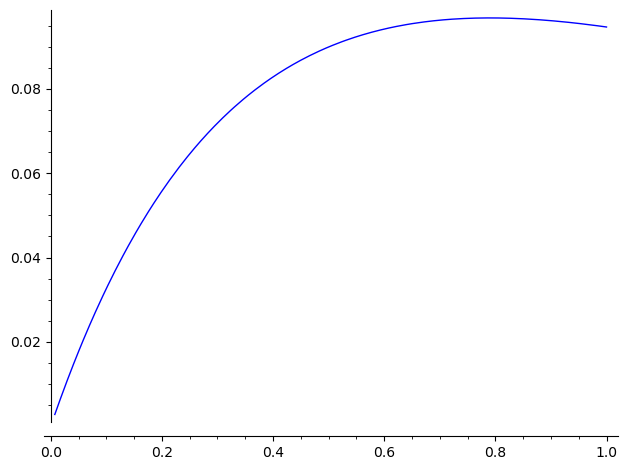}$\qquad$
\includegraphics[width=0.40\textwidth]{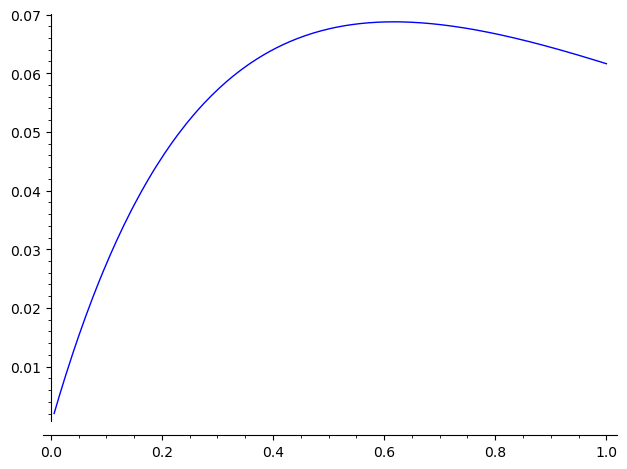}\\
\includegraphics[width=0.40\textwidth]{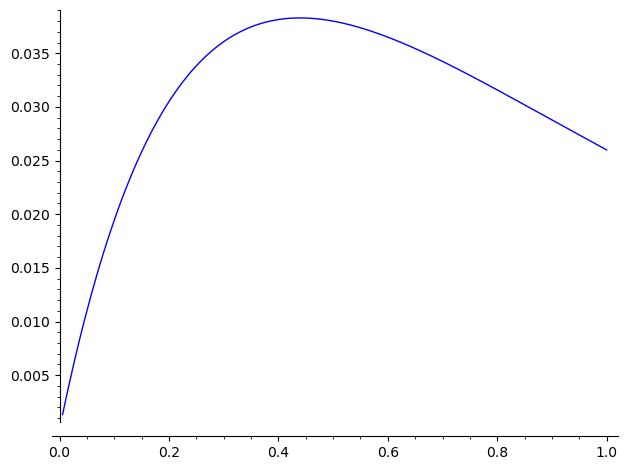}$\qquad$
\includegraphics[width=0.40\textwidth]{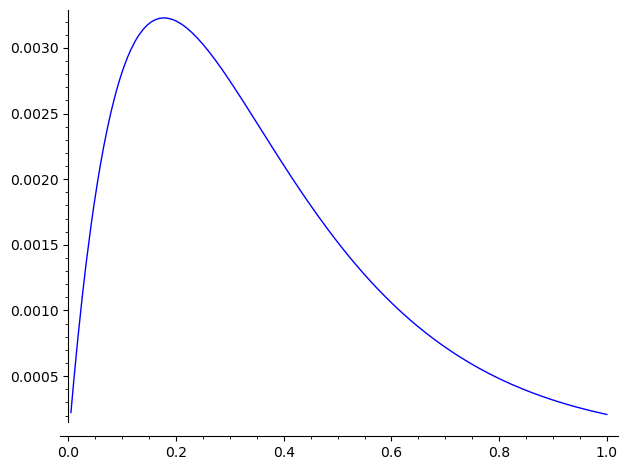}\\
\caption{\sl\footnotesize Plot of~$\left(0,1\right)\ni s\mapsto
{\mathcal{G}}_2(s;L,T)$ (say, when~$T=1$)
for~$L\in\{1,1.2,1.3,1.5,2,10 \}$.}
        \label{63wDIK-MD8990dfgh35dede7-04yeufLO}
\end{figure}
\end{center}

The functional~${\mathcal{G}}_2$ also shows
an interesting bifurcation diagram plotted in Figure~\ref{63wDIK-MD8990dfgh35dede7-04yeufLO}:
also in this framework preys located close to the origin favor
local diffusive strategies (optimized in this case for~$s=1$)
and when~$L$ becomes larger an larger the optimal exponent
moves to the left till becomes~$s=0$. Analogously to the functional~${\mathcal{G}}_1(s;L,T)$, we can show that
\begin{equation}\label{SLL-1.K}
\begin{split}&
{\mbox{if $L$ is large enough, the optimal foraging strategy
for~${\mathcal{G}}_2(s;L,T)$}}\\ &{\mbox{is uniquely attained at some~$s_L$ such that }}\lim_{L\to+\infty}
s_L=0.\end{split}
\end{equation}

The claim in \eqref{SLL-1} is proved analytically in Appendix \ref{ALLc-1.kk}. Note that also in this case we have the same instability of the optimal search strategy that we already observed for $\mathcal{G}_1$.

The similarities and differences between
Figures~\ref{63wDIK-MD89907-04yeufLO} and~\ref{63wDIK-MD8990dfgh35dede7-04yeufLO}
highlight how different normalization
choices in the model can affect optimal strategies: note indeed
that the only difference
between~${\mathcal{G}}_1$ and~${\mathcal{G}}_2$ lies
in the way the diffusion coefficient~$\kappa$  
is modeled on the basis of the underlying random process.
The sensitivity of the optimization strategies
with respect to these normalizing constants seems to be not investigated
in the current literature
and it produces in
Figures~\ref{63wDIK-MD89907-04yeufLO}
and~\ref{63wDIK-MD8990dfgh35dede7-04yeufLO}
a different outcome on the optimality of the Gaussian exponent~$s=1$;
this interesting difference
is induced by the analytical observation that~${\mathcal{G}}_1(1;L,T)=0<
{\mathcal{G}}_2(1;L,T)$.

\begin{center}
\begin{figure}[!ht]
\includegraphics[width=0.40\textwidth]{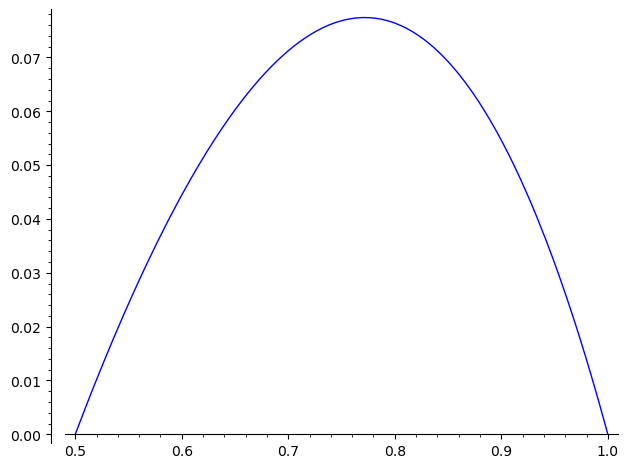}$\qquad$
\includegraphics[width=0.40\textwidth]{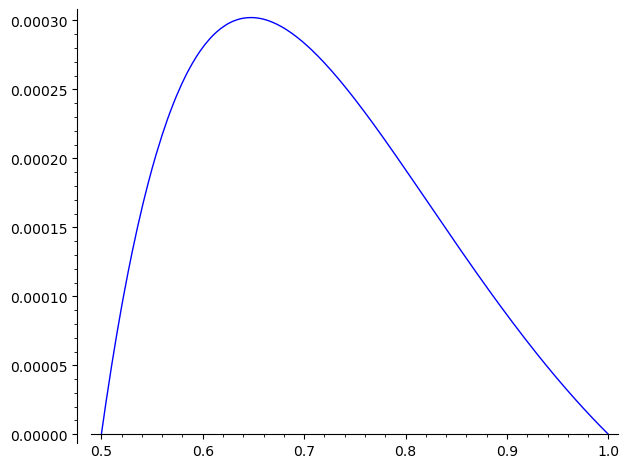}\\
\includegraphics[width=0.40\textwidth]{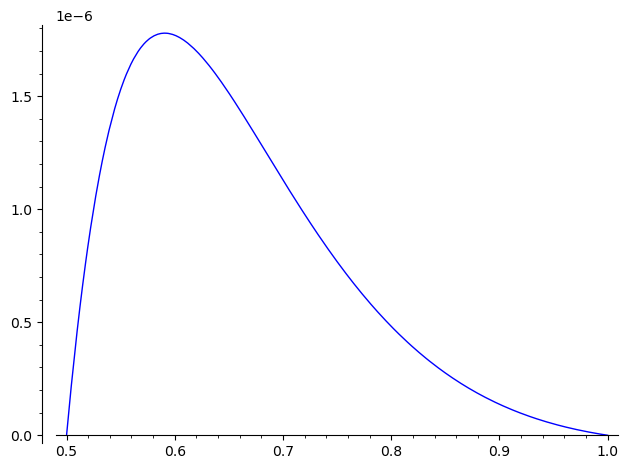}$\qquad$
\includegraphics[width=0.40\textwidth]{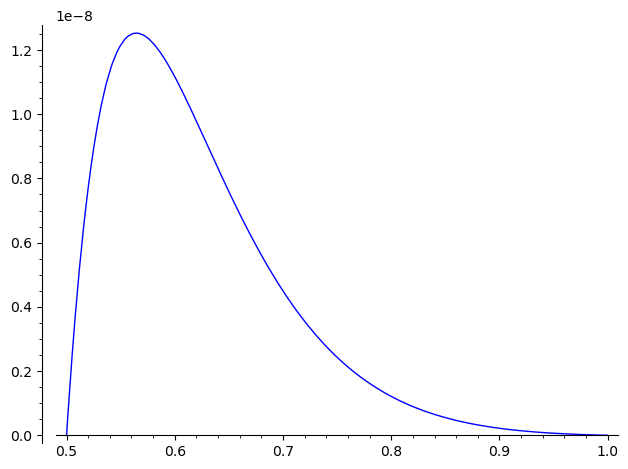}\\
\includegraphics[width=0.40\textwidth]{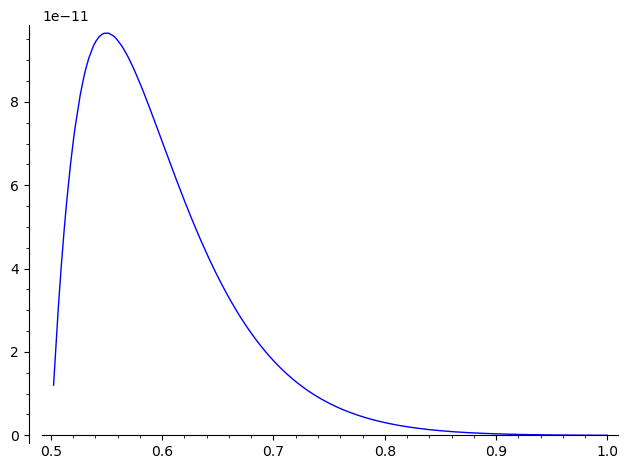}$\qquad$
\includegraphics[width=0.40\textwidth]{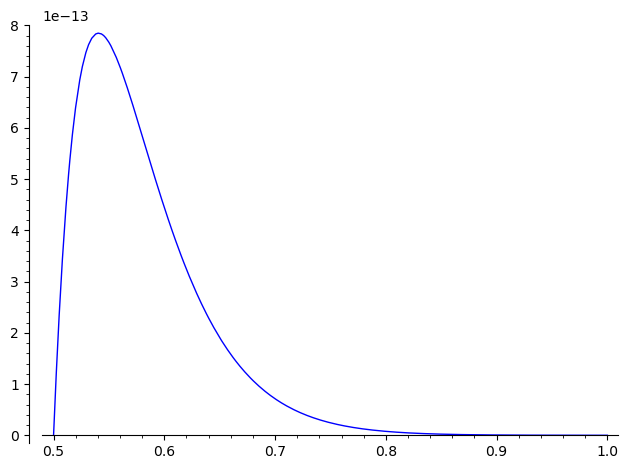}
\caption{\sl\footnotesize Plot of~$\left(\frac12,1\right)\ni s\mapsto
{\mathcal{G}}_3(s;L,T)$ for $T=1$ and~$L=10^j$, $j\in\{0,\dots,5\}$.}
        \label{6tgh3wDIK-MDye09876543nijufLO}
\end{figure}
\end{center}

As for the functional~${\mathcal{G}}_3$, plots for
different values of~$L$ are given in Figure~\ref{6tgh3wDIK-MDye09876543nijufLO}.
Interestingly, on the one hand,
both the inverse square law~$s=\frac12$
and the Gaussian law~$s=1$ are minima for the functional for every~$T$ and~$L$;
on the other hand, for very sparse targets (corresponding to large values of~$L$),
\begin{equation}\label{SLL}
\begin{split}&
{\mbox{the optimal foraging strategy
for~${\mathcal{G}}_3(s;L,T)$ is uniquely attained at some~$s_L$}}\\
&{\mbox{such that }}\lim_{L\to+\infty}
s_L=\frac12.\end{split}\end{equation}
This can also be proved analytically, see Appendix~\ref{ALL}.

\begin{center}
\begin{figure}[!ht]
\includegraphics[width=0.40\textwidth]{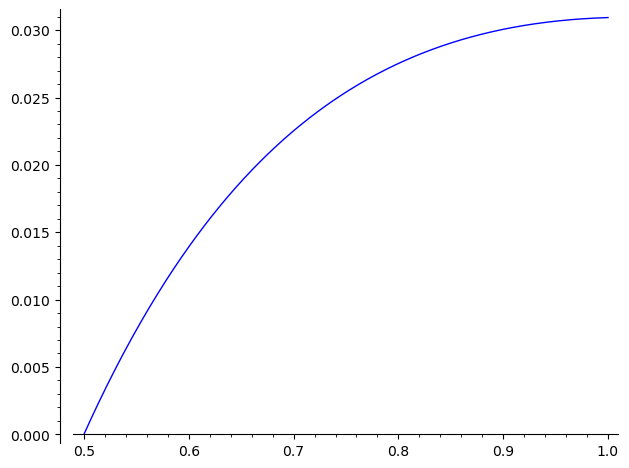}$\qquad$
\includegraphics[width=0.40\textwidth]{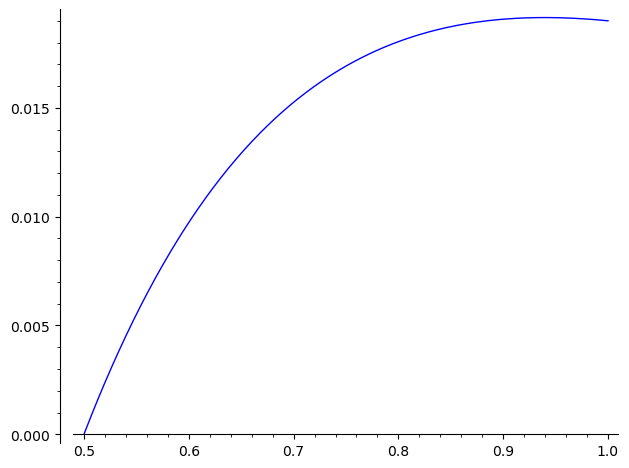}\\
\includegraphics[width=0.40\textwidth]{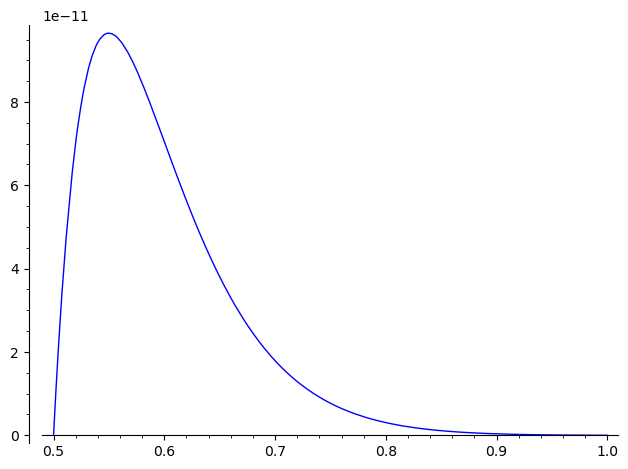}$\qquad$
\includegraphics[width=0.40\textwidth]{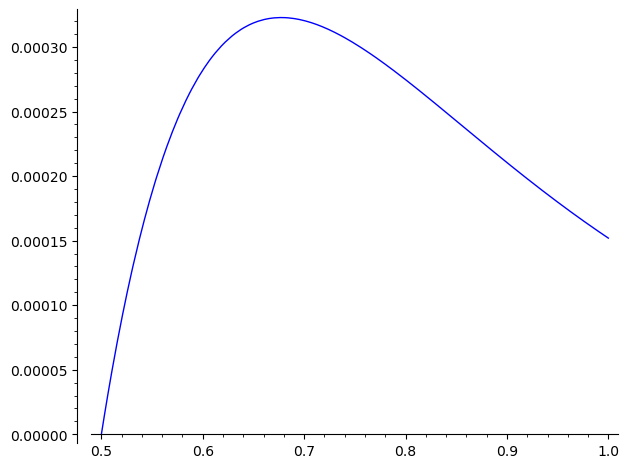}\\
\includegraphics[width=0.40\textwidth]{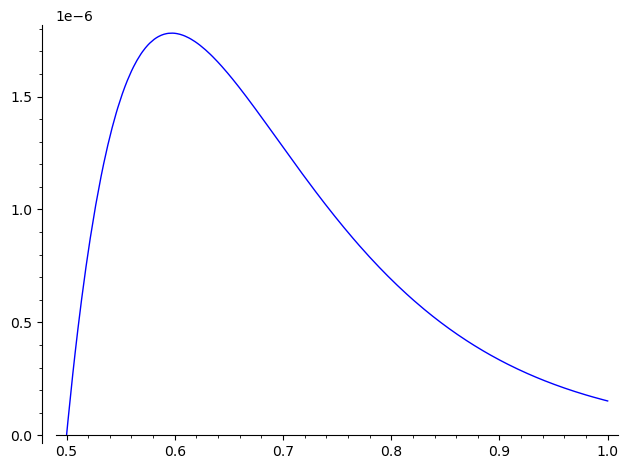}$\qquad$
\includegraphics[width=0.40\textwidth]{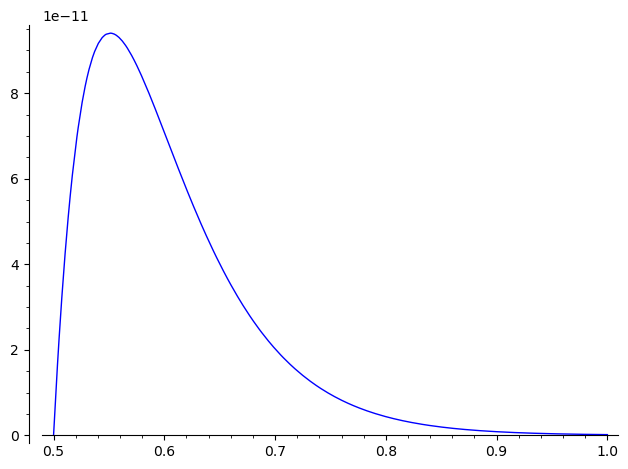}
\caption{\sl\footnotesize Plot of~$\left(\frac12,1\right)\ni s\mapsto
{\mathcal{G}}_4(s;L,T)$ for $T=1$ and~$L\in\{1.7,2,4,10,100,10000\}$.}
        \label{6tgh3wDIK-MDye09876543nijufLO-06557}
\end{figure}
\end{center}

The graph of the functional~${\mathcal{G}}_4$
is instead plotted in Figure~\ref{6tgh3wDIK-MDye09876543nijufLO-06557}:
notice that for~$L\leq 1.7$ this functional is increasing
and attains its maximum for the Gaussian strategy~$s=1$,
but for ~$L\geq 2$ the functional~${\mathcal{G}}_4$ develops an interior maximum. More precisely, as proved analytically in Appendix~\ref{A4444}, setting
\begin{equation}\label{ELLSTARVO} L^\star:=\exp\left(-\frac{\zeta'(2)}{\zeta(2)}\right)= 1.768198...\end{equation}
we have that     
\begin{equation}\label{ASOG4}
\begin{split}
&{\mbox{if~$L\leq L^\star$ the optimal foraging strategy
for~${\mathcal{G}}_4(s;L,T)$ is uniquely attained at~$s=1$,}}\\
&{\mbox{while if~$L>L^\star$ the optimal foraging strategy
for~${\mathcal{G}}_4(s;L,T)$ is uniquely attained}}\\
&{\mbox{at some }} s_L\in \left(\frac12,1\right) {\mbox{
such that}} \lim_{L\to+\infty}s_L=\frac12.
\end{split}
\end{equation}
Therefore, for larger and larger values of $L>L^\star$, the optimal strategy for $\mathcal{G}_4$ is getting closer and closer to the inverse power law distribution induced by $s=\frac12$.

\begin{center}
\begin{figure}[h]
\includegraphics[width=0.40\textwidth]{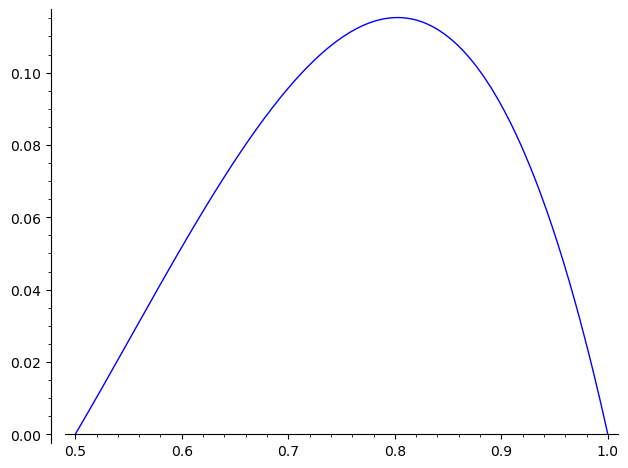}
\caption{\sl\footnotesize Plot of~$\left(\frac12,1\right)\ni s\mapsto
{\mathcal{G}}_5(s;L,T)$ for $
L=T^{\frac{2s-1}{2s(1+2s)}}$.}
        \label{6tgSQRTGHS78dh3wDIK-MDye09876543nijufLO-06557}
\end{figure}
\end{center}

\begin{center}
\begin{figure}[h]
\includegraphics[width=0.40\textwidth]{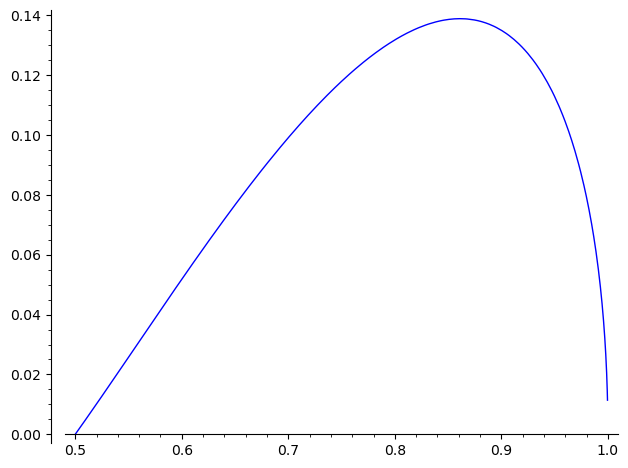}
\caption{\sl\footnotesize Plot of~$\left(\frac12,1\right)\ni s\mapsto
{\mathcal{G}}_6(s;L,T)$ for $
L=T^{\frac{2s-1}{2s(1+2s)}}$.}
        \label{6tgSQRTG567890hhh-HS78dh3wDIK-MDye09876543nijufLO-06557}
\end{figure}
\end{center}

The cases of~$
{\mathcal{G}}_5$ and~$
{\mathcal{G}}_6$
are quite sophisticated, since their optimization strategies depend both on the final time~$T$
and on the scantness of the targets modeled by~$L$.
In the special situation in which~$
L=T^{\frac{2s-1}{2s(1+2s)}}$ these value functionals do not depend on~$L$
and~$T$ and they are plotted in Figures~\ref{6tgSQRTGHS78dh3wDIK-MDye09876543nijufLO-06557}
and~\ref{6tgSQRTG567890hhh-HS78dh3wDIK-MDye09876543nijufLO-06557}, respectively.

In particular, it appears that when~$
L=T^{\frac{2s-1}{2s(1+2s)}}$ the functional~$
{\mathcal{G}}_5$
is maximized in proximity of~$s=
0.80261...$ and~${\mathcal{G}}_6$
in proximity of~$s=
0.861187...$: in spite of the rather arbitrary choice
relating~$T$ and~$L$ in
Figures~\ref{6tgSQRTGHS78dh3wDIK-MDye09876543nijufLO-06557}
and~\ref{6tgSQRTG567890hhh-HS78dh3wDIK-MDye09876543nijufLO-06557},
it is suggestive to compare these optimal intermediate values
between the inverse square power law distribution and the Gaussian one
with the ones observed experimentally
for some blue sharks
(see e.g. Figure~1b in~\cite{HUM}
which would correspond to~$s=0.73$),
basking sharks
and bigeye tunas
(see e.g. Figures~1b and 1c in~\cite{SCA}
which would correspond to~$s=0.7$). 
See also the red curves in Figures~2, 4
and~6 of~\cite{Palyulin2931}
(which corresponds to~$s\simeq0.75$).
It is also interesting to compare these values with the
simulation data of some swarm dynamics
(corresponding to~$s=0.745$, see Figure~3
in~\cite{SWA}).
Of course, we are not aiming here at precisely reconstructing
the quantitative results arising in specific real-world experiments,
but we think it is an interesting feature that
even the very simplified situation that we describe may
lead to optimal values of~$s$ which are
somewhat intermediate between~$s=\frac12$
and~$s=1$.

\subsection{Prey at the origin and remote prey}

We now consider the case of two targets, one located ``at home''
at the origin and another far away at a
given distance~$L>0$. This corresponds to a prey distribution of the form
$$ p_{0,L}(x):=\delta_0(x)+\delta_L(x).$$
Since the foraging success functional in~\eqref{AV}
is linear with respect to the target distribution,
the analysis of this case reduces to the superposition
of the value functionals introduced in~\eqref{UTIFU}
and~\eqref{LFT}: thus, in the above notation, we define
\begin{equation}\label{MK:03} {\mathcal{H}}_j(s;L,T):=
{\mathcal{E}}_j(s;T)+{\mathcal{G}}_j(s;L,T)
\qquad{\mbox{for }}j\in\{1,\dots,6\}\end{equation}
and we find that
\begin{equation*}
\begin{split}&{\mathcal{H}}_1(s;L,T)=\frac{1}{\pi \,T^{\frac{1}{2s}}\,(2s-1)}
\Gamma\left( \frac{1}{2s}\right)+
\frac{T \;\Gamma(1+2s)\,\sin(\pi s) }{{2}\pi \,L^{1+2s}},\\
&
{\mathcal{H}}_2(s;L,T)=\frac{2 \,\big(-2\zeta(-2s)\big)^{\frac{1}{2s}}}{T^{\frac{1}{2s}}\,(2s-1)}
\Gamma\left( \frac{1}{2s}\right)+
\frac{T}{{4}L^{1+2s}\, \zeta(1+2 s)}
 ,\\
&
{\mathcal{H}}_3(s;L,T)=
\frac{\zeta(1+2s)}{\pi \,T^{\frac{1}{2s}}\,(2s-1)\,\zeta(2s)}
\Gamma\left( \frac{1}{2s}\right)+
\frac{T \;\zeta(1+2s)\,\Gamma(1+2s)\,\sin(\pi s) }{{2}\pi L^{1+2s}\,\zeta(2s)},\\
&{\mathcal{H}}_4(s;L,T)=
\frac{\zeta(1+2s)2 \big(-2\zeta(-2s)\big)^{\frac{1}{2s}}
}{
 T^{\frac{1}{2s}}\,(2s-1)\,
\zeta(2s)}
\Gamma\left( \frac{1}{2s}\right)+
\frac{T }{{4} L^{1+2s}\zeta(2 s)},\\
& {\mathcal{H}}_5(s;L,T)
=
\frac{(1+2s)\;\Gamma\left( \frac{1}{2s}\right)}{4 s \,(2s-1)\,T^{\frac{1}{s}}\;
\Gamma\left(\frac{2s-1}{2s}\right)}+
\frac{T^{\frac{2s-1}{2s}} \;\Gamma(2+2s)\,\sin(\pi s) }{{8} L^{1+2s}\,\Gamma\left(\frac{2s-1}{2s}\right)}
\\{\mbox{and }}\qquad
&
{\mathcal{H}}_6(s;L,T)=
\frac{\pi^2\big(-2\zeta(-2s)\big)^{\frac{1}{s}}\,(1+2s)\;\Gamma\left( \frac{1}{2s}\right)}{ s \,(2s-1)\,T^{\frac{1}{s}}\;
\Gamma\left(\frac{2s-1}{2s}\right)}\\&\qquad\qquad\qquad\qquad
+
\left( \frac{-1}{2^{1+2s}\pi^{2s}\,\zeta(-2s)}\right)^{\frac{2s-1}{2s}}
\frac{T^{\frac{2s-1}{2s}} \;\Gamma(2+2s)\,\sin(\pi s) }{{8} L^{1+2s}\,
\Gamma\left(\frac{2s-1}{2s}\right)}.
\end{split}\end{equation*}
We point out that in all the above value functionals, the second term
becomes predominant for large values of~$T$, hence
the long time analysis for~${\mathcal{H}}_j$
boils down to the one developed for~${\mathcal{G}}_j$
in Section~\ref{REMOTE} (this is consistent with the idea that
for long times the forager has drifted away from its initial location).
Similarly, small values of~$T$ reduce the analysis 
of~${\mathcal{H}}_j$
to the one developed for~${\mathcal{E}}_j$
in Section~\ref{SINOR} (consistently with the ansatz that for small
times the forager will exploit the targets in the vicinity of its burrow).

\begin{center}
\begin{figure}[!ht]
\includegraphics[width=0.40\textwidth]{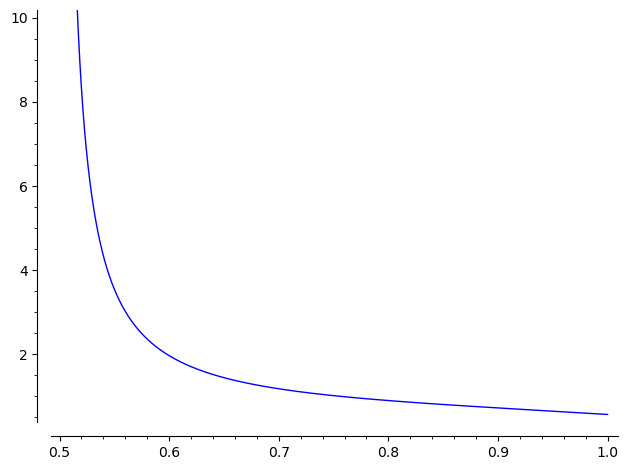}$\qquad$
\includegraphics[width=0.40\textwidth]{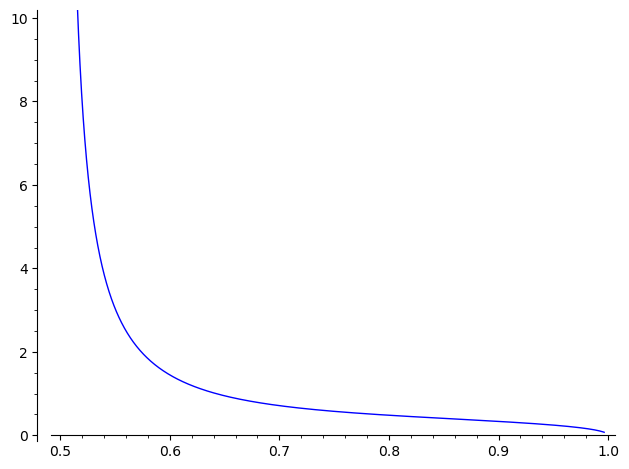}\\
\includegraphics[width=0.40\textwidth]{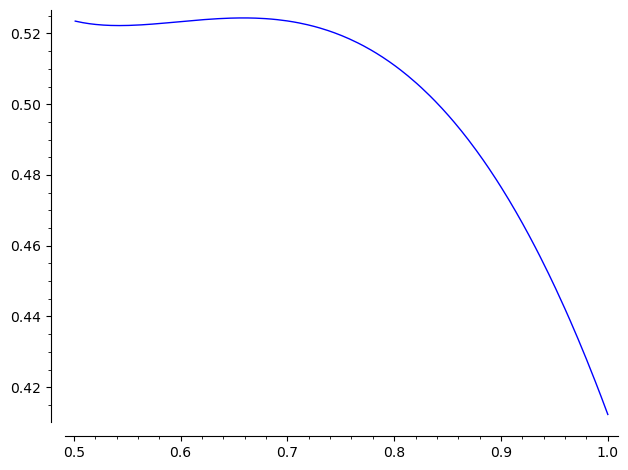}$\qquad$
\includegraphics[width=0.40\textwidth]{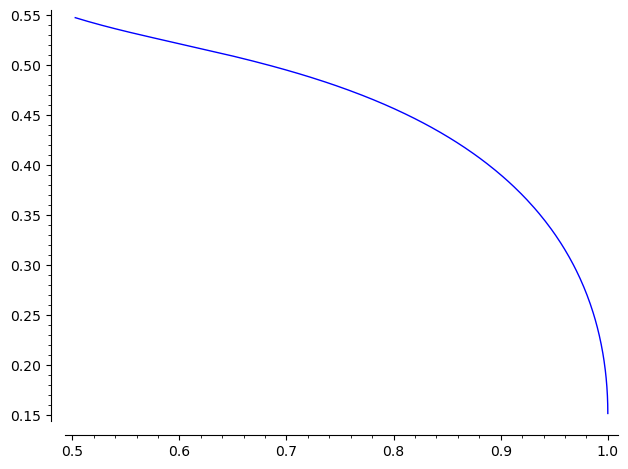}\\
\includegraphics[width=0.40\textwidth]{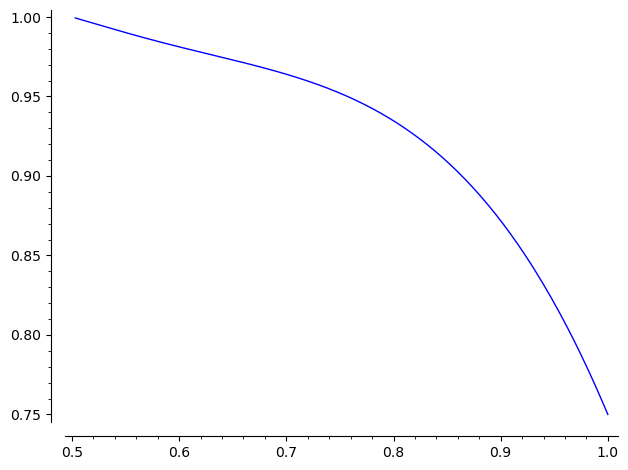}$\qquad$
\includegraphics[width=0.40\textwidth]{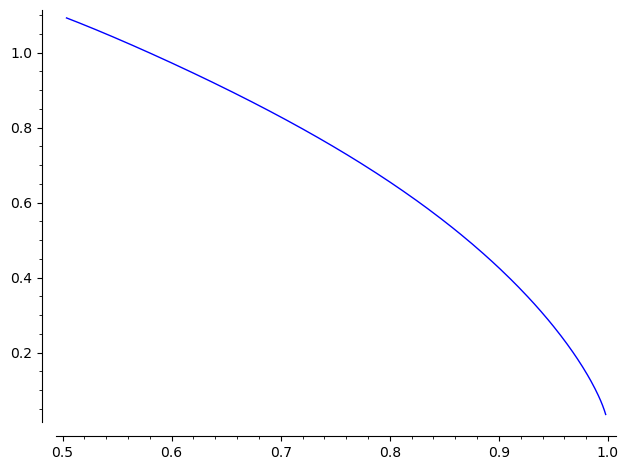}
\caption{\sl\footnotesize Plot of~$\left(\frac12,1\right)\ni s\mapsto
{\mathcal{H}}_j(s;L,T)$ for $T=1=L$ and~$j\in\{1,\dots,6\}$.}
        \label{6tgh3wDIK-MDye0679-08870984ndjgfg-133-49876543nijufLO-06557}
\end{figure}
\end{center}

Instead, when~$T=1$ both the terms in~${\mathcal{H}}_j$
contribute to the optimization of~${\mathcal{H}}_j$ and the corresponding plots (when also~$L=1$) are given in Figure~\ref{6tgh3wDIK-MDye0679-08870984ndjgfg-133-49876543nijufLO-06557}, showing an optimal
foraging strategy corresponding to~$s=\frac12$ in these specific situations.

\section*{Acknowledgments}

SD and EV are members of AustMS.
Supported by the Australian Laureate Fellowship FL190100081 ``Minimal
surfaces, free boundaries and partial differential equations''
and
by the Australian Research Council DECRA DE180100957
``PDEs, free boundaries
and applications''.

\begin{appendix}

\section{Appendices}

\subsection{Analytical verification of~\eqref{UNST}}\label{UNSTAPP}

Let
\begin{equation}\label{CIESTD} c(s;T):=\frac{1}{2\pi T^{\frac{1}{2s}}(2s-1)^2 s^2}\,\Gamma\left(\frac{1}{2s}\right)\end{equation}
and notice that
\begin{equation}\label{C55}
\begin{split}
\frac{d}{ds}c(s;T)\,&=\,\left(\frac{\ln T}{2s^2}-\frac{4}{2s-1}-\frac2s\right)\,c(s;T)
-\frac{1}{4\pi T^{\frac{1}{2s}}(2s-1)^2 s^4}\,\Gamma'\left(\frac{1}{2s}\right)
\\&=\,\left(\frac{\ln T}{2s^2}-\frac{4}{2s-1}-\frac2s-\frac{1}{2s^2}\psi\left(\frac1{2s}\right)\right)\,c(s;T).
\end{split}\end{equation}
Moreover, by~\eqref{UTIFU},
\begin{equation}\label{ZIO}
{\mathcal{E}}_1(s;T)=
2s^2(2s-1)\,c(s;T)
\end{equation}
and therefore
\begin{equation}\label{eq_1}\begin{split}
\frac{d}{ds}\mathcal{E}_1(s;T)
\,&=\,\big(4s (2s-1)+4s^2\big)\,c(s;T)+
2s^2(2s-1)
\frac{d}{ds}c(s;T)
\\&=\,
\left((2s-1) \ln T-4s^2
-(2s-1)\psi\left(\frac1{2s}\right)\right)\,c(s;T)
.\end{split}\end{equation}
Now we define
\begin{equation}\label{ZETADIX}
Z(x):=\sum_{k=1}^{\infty}\left(\frac{1}{k+x}-\frac{1}{k}\right) \end{equation}
and (see Chapter~10 in~\cite{GAM}) we find that
\begin{equation}\label{RViNSXC}
\psi(x)=-\gamma-\frac1x-\sum_{k=1}^{+\infty}\left(
\frac{1}{k+x}-\frac1k
\right)=-\gamma-\frac1x-Z(x).
\end{equation}
Accordingly, we can rewrite equation \eqref{eq_1} as
\begin{equation}\label{ZSTI0}
\frac{d}{ds}\mathcal{E}_1(s;T)=
\left((2s-1) \ln T
+(2s-1)\gamma -2s+(2s-1)\,Z\left(\frac1{2s}\right)
\right)\,c(s;T).
\end{equation}
Moreover, since~$s\in \left(\frac{1}{2},1\right)$, we have that~$\frac1{2s}\in\left(\frac{1}{2},1\right)$
and therefore, for each~$k\ge1$,
$$ \frac{1}{k+\frac1{2s}}-\frac1k\ge\frac{1}{k+1}-\frac1k=-\frac{1}{k(k+1)}\ge-\frac{1}{k^2}.
$$
As a result,
\begin{equation}\label{ZSTI}
-1<-\frac{\pi^2}{12}=-\frac{1}{2}\sum_{k=1}^{\infty}\frac{1}{k^2}\leq \frac12 Z\left(\frac{1}{2s}\right)\leq 0.
\end{equation}
We now define the function
\begin{equation}\label{fanbdg}
f(s):=
\Big((2s-1) \ln T
+(2s-1)\gamma-2s -2(2s-1)
\Big)\,c(s;T)
\end{equation}
and we deduce from~\eqref{ZSTI0}
and~\eqref{ZSTI} that
\begin{equation}\label{TRAm}
f(s)\leq \frac{d}{ds}\mathcal{E}_1(s;T)
.\end{equation}
For~$T$ large enough such that~$\ln T+\gamma\ge3$, let
\begin{equation}\label{PDA0}
\overline{s}_T:=\frac{1}{2}\left(\frac{\ln T+\gamma-2}{\ln T+\gamma-3}\right)
\end{equation}
and observe that 
\begin{equation}\label{PDA}
f(s)\geq 0 {\mbox{ if and only if }} s \in \left(\overline{s}_T,1\right).
\end{equation}
Recalling \eqref{UTIFU}, we observe that 
$$
\lim_{s\rightarrow \frac{1}{2}}\mathcal{E}_1(s,T)=+\infty.
$$
This limit together with \eqref{PDA} and~\eqref{TRAm} imply that
\begin{equation}\label{DOVED}
{\mbox{the map $\left(\frac12,1\right)\ni s\mapsto\mathcal{E}_1(s;T)$ attains its minimum 
somewhere in the interval $\bigg(\frac{1}{2},\overline{s}_T\bigg)$.}}\end{equation}
Furthermore, by~\eqref{C55} and~\eqref{eq_1},
\begin{eqnarray*}&&
\frac{d^2}{ds^2}\mathcal{E}_1(s;T)\\
&=&\left(2\ln T-8s
-2\psi\left(\frac1{2s}\right)+\frac{2s-1}{2s^2}\psi'\left(\frac1{2s}\right)\right)\,c(s;T)\\&&
+
\left((2s-1) \ln T-4s^2
-(2s-1)\psi\left(\frac1{2s}\right)\right)\,\frac{d}{ds}c(s;T)\\
&=&
\Bigg[\left(2\ln T-8s
-2\psi\left(\frac1{2s}\right)+\frac{2s-1}{2s^2}\psi'\left(\frac1{2s}\right)\right)\\&&
+
\left((2s-1) \ln T-4s^2
-(2s-1)\psi\left(\frac1{2s}\right)\right)\left(\frac{\ln T}{2s^2}-\frac{4}{2s-1}-\frac2s-\frac{1}{2s^2}\psi\left(\frac1{2s}\right)\right)\Bigg]c(s;T)\\
&\ge&
\Bigg[
\frac{(2s-1)\ln^2T}{2 s^2}-
\left(
4+\frac{2(2s-1)}s+\frac{2s-1}{s^2}\psi\left(\frac1{2s}\right)
\right)\ln T
+\frac{16 s^2}{2 s - 1}-C
\Bigg]c(s;T)
\end{eqnarray*}
for some~$C>0$ independent on~$T$.

Hence, 
setting~$A:=\frac{\sqrt{2s-1}\,\ln T}{\sqrt{2}\, s}$ and~$B:=\frac{4s}{\sqrt{2 s - 1}}$, and
noticing that
$$
4\sqrt{2}\,\ln T
=2AB\le A^2+B^2=
\frac{(2s-1)\ln^2T}{2 s^2}+\frac{16 s^2}{2 s - 1},$$
we conclude that
\begin{equation}\label{MpsodkKSA}
\frac{d^2}{ds^2}\mathcal{E}_1(s;T)
\ge
\Bigg[
\left(4\sqrt2-
4-\frac{2(2s-1)}s-\frac{2s-1}{s^2}\psi\left(\frac1{2s}\right)
\right)\ln T -C
\Bigg]c(s;T)
.\end{equation}
Now, recalling~\eqref{RViNSXC}
and~\eqref{ZSTI},
\begin{eqnarray*}
&&4\sqrt2-
4-\frac{2(2s-1)}s-\frac{2s-1}{s^2}\psi\left(\frac1{2s}\right)\\
&=&4\sqrt2-
4-\frac{2(2s-1)}s+\frac{2s-1}{s^2}\left(
\gamma+2s+Z\left(\frac1{2s}\right)
\right)\\&=&4\sqrt2-
4+\frac{2s-1}{s^2}\left(
\gamma+Z\left(\frac1{2s}\right)
\right)\\&\ge&4\sqrt2-
4-\frac{2s-1}{s^2}\left(2-\gamma\right)\\&\ge&4\sqrt2-
4-2+\gamma,
\end{eqnarray*}
which is strictly positive. 

This observation and~\eqref{MpsodkKSA} entail that~$\mathcal{E}_1$
is strictly convex for large~$T$, and
consequently, the minimum in~\eqref{DOVED}
is attained at a unique point that we denote by~$s_T$.
Since
$$
\lim_{T\rightarrow\infty}\overline{s}_T=\frac{1}{2},$$
we thus obtain from~\eqref{DOVED} that
$$\lim_{T\rightarrow\infty}s_T=\frac{1}{2},
$$
which concludes the analytical proof of~\eqref{UNST}.

\subsection{Analytical verification of~\eqref{UNST-2}}\label{VF}

We point out that
\begin{equation}\label{LimitE_2}
\lim_{s\searrow1/2}{\mathcal{E}}_2(s;T)=+\infty,
\end{equation}
since the denominator in the definition of~${\mathcal{E}}_2$
in~\eqref{UTIFU} vanishes. Moreover,
since the Riemann Zeta Function
vanishes at the negative even integers (different from~$-1$), we note that $\zeta(-2)=0$, whence
\begin{equation}\label{ESPOIHGFSVBSU1213}
\lim_{s\nearrow1}{\mathcal{E}}_2(s;T)=0.\end{equation}
Also,
in light of~\eqref{UTIFU}, we can rewrite ${\mathcal{E}}_2(s;T)$ as
\begin{equation*}
{\mathcal{E}}_2(s;T)=2\pi (h(s))^{\frac{1}{2s}}\mathcal{E}_1(s,T),
\end{equation*}
where we have defined $h(s):=-2\zeta(-2s)>0$ for each $s\in (1/2,1)$.

We remark\footnote{Here we are using the Riemann's Functional Equation
$$ \zeta(z)=2^z\pi^{z-1}\sin\left(\frac{\pi z}2\right)\Gamma(1-z)\zeta(1-z).
$$} that, if~$s\in (1/2,1)$,
\begin{eqnarray*}
h(s)=2^{1-2s}\pi ^{-2s-1}\, \sin (\pi s)\ \Gamma (1+2s)\, \zeta (1+2s)\in
\left(0,\frac{ \Gamma (3)\, \zeta (2)}{\pi ^2}\right]=\left(0,\frac{1}{3}\right]
\end{eqnarray*}
and therefore
\begin{equation}\label{EMEELOD} \ln(h(s))\le \ln\frac13<0.\end{equation}
Also,
\begin{equation}\label{EGGIU1}
\begin{split}
h'(s)\,&=\,
2^{1 - 2 s} \pi^{-2 s - 1} \Gamma(1+2 s) \Big[ 2 \sin(\pi s) \zeta'(1+2 s) \\&\qquad+ \zeta(1+2 s) \Big(
\pi \cos(\pi s) + 2 \sin(\pi s)\big( \psi(1+2 s)-\ln2 - \ln\pi\big)\Big)\Big].
\end{split}
\end{equation}
We recall that
$$ \psi(1+2 s)-\ln2 - \ln\pi\le\psi(3)-\ln2 - \ln\pi=
-0.9150927...$$
and therefore we infer from~\eqref{EGGIU1}
that
\begin{equation*}
\begin{split}
h'(s)\,&\le\,
2^{1 - 2 s} \pi^{-2 s - 1} \Gamma(1+2 s ) \Big[ 2 \sin(\pi s) \zeta'(1+2 s)+ \pi \zeta(1+2 s)\cos(\pi s)
\Big]\\&\le\,
2^{1 - 2 s} \pi^{-2 s } \Gamma(1+2 s)  \zeta(1+2 s)\cos(\pi s).
\end{split}
\end{equation*}
Since the latter term is nonpositive for~$s\in\left[\frac12,1\right]$, and actually 
strictly negative when~$s\in\left(\frac12,1\right]$, the observations above yield that~$h'<0$
in~$\left(\frac12,1\right)$ and more precisely, for all~$ s\in \left(\frac{1}{2},1\right)$,
\begin{equation}\label{YOPH}
-C_1<h'(s)<-C_2 ,
\end{equation}
where $C_1$ and~$C_2$ are positive constants.

We now compute the derivative with respect to $s$ of ${\mathcal{E}}_2(s;T)$ and get
\begin{equation*}
\frac{d}{ds}{\mathcal{E}}_2(s;T)=2\pi\left[\frac{d}{ds}\mathcal{E}_1(s,T)(h(s))^{\frac{1}{2s}}
+\mathcal{E}_1(s,T)(h(s))^{\frac{1}{2s}}\left(-\frac{1}{2s^2}\ln(h(s))+\frac{1}{2s}\frac{h'(s)}{h(s)}       \right)\right].
\end{equation*}
Hence, from equations~\eqref{ZIO} and~\eqref{ZSTI0} we deduce that
\begin{equation}\label{A345NS-03e30}
\frac{d}{ds}{\mathcal{E}}_2(s;T)=2\pi c(s;T)(h(s))^{\frac{1-2s}{2s}}\Big(
P(s,T)\,h(s)-(2s-1)\,\ln(h(s))\,h(s)+s(2s-1)\,h'(s)\Big),
\end{equation}  
where
$$
P(s,T):=(2s-1) \ln T
+(2s-1)\gamma -2s+(2s-1)\,Z\left(\frac1{2s}\right)
$$
and~$c$ and~$Z$ are as defined in~\eqref{CIESTD} and~\eqref{ZETADIX}, respectively.

Now, using~\eqref{ZSTI} and~\eqref{fanbdg}, we observe that
\begin{eqnarray*}
P(s,T)&=&(2s-1) \ln T
+(2s-1)\gamma -2s+(2s-1)\,Z\left(\frac1{2s}\right)\\& =&
(2s-1) \ln T
+(2s-1)\gamma -2s
-2(2s-1) +(2s-1)\left(2+Z\left(\frac1{2s}\right)\right)
\\&\ge& (2s-1) \ln T
+(2s-1)\gamma -2s
-2(2s-1)\\&=&\frac{f(s)}{c(s;T)}.
\end{eqnarray*}
As a consequence, owing also to~\eqref{EMEELOD},~\eqref{YOPH} and~\eqref{A345NS-03e30}
we find that
\begin{equation}\label{Pop}\begin{split}
\frac{d}{ds}{\mathcal{E}}_2(s;T)\ge\,&
2\pi c(s;T)\,(h(s))^\frac{1-2s}{2s}
\left(\frac{f(s)\,h(s)}{c(s;T)}+s(2s-1)\,h'(s)
\right)\\ \ge\,&
2\pi c(s;T)\,(h(s))^\frac{1-2s}{2s}\left(\frac{f(s)\,h(s)}{c(s;T)}-s(2s-1)
C_1 \right).
\end{split}\end{equation}
Consequently, by~\eqref{PDA} and the fact that $h(s)$ is decreasing,
we get that
for all~$ s \in \left(\overline{s}_T,\frac{3}{4}\right)$,
\begin{equation}\label{pizza5}
2\pi c(s;T)h(s)^\frac{1-2s}{2s}\left(\frac{f(s)\,h(3/4)}{c(s;T)}-s(2s-1)C_1\right)\leq \frac{d}{ds}{\mathcal{E}}_2(s;T) 
\end{equation}
where $\overline{s}_T$ is as defined in~\eqref{PDA0}. Now we set 
\begin{equation*}
\hat{s}_T:= \overline{s}_T+\frac{1}{2}\left(\frac{s(2s-1)
C_1}{h(3/4)\big(\ln T-3+\gamma  \big)}\right)=
\frac{1}{2}\left(\frac{\ln T-2+\gamma+s(2s-1)
C_1/h(3/4)}{\ln T-3+\gamma}\right).
\end{equation*} 
Hence, recalling again~\eqref{fanbdg} and assuming~$T$ conveniently large, 
for all~$
s \in \left(\hat{s}_T,\frac{3}{4}\right)$,
\begin{equation*}
\begin{split}
\frac{f(s)h(3/4)}{c(s;T)}-s(2s-1)C_1\,&=\,
h(3/4)\Big((2s-1) \ln T
+(2s-1)\gamma-2s -2(2s-1)\Big)
-s(2s-1)C_1\\&=\,
2s\,h(3/4)\Big( \ln T
+\gamma-3
\Big)
+ h(3/4)\Big(-\ln T-\gamma+2\Big)-s(2s-1)
C_1\\
&>\,2\hat{s}_T\,h(3/4)\Big( \ln T
+\gamma-3\Big)
+h(3/4)\Big(- \ln T-\gamma+2\Big)-s(2s-1)C_1\\&=\,0.
\end{split}
\end{equation*}
Plugging this information into~\eqref{pizza5}
we find that, for large~$T$, the function~$\left(\frac{1}{2},1\right)\ni s\mapsto
{\mathcal{E}}_2(s;T)$
is strictly increasing in~$\left(\hat{s}_T,\frac{3}{4}\right)$ and
in view of \eqref{LimitE_2} we obtain that this function possesses a local minimum~$s_T\in
\big(\frac{1}{2},\hat{s}_T\big]$ and (recalling also~\eqref{ESPOIHGFSVBSU1213})
a local maximum~$S_T\in\left[\frac34,1\right)$.
Since $\hat{s}_T$  converges to~$\frac{1}{2}$ as~$T\to+\infty$, we also
obtain that~$s_T\to\frac12$ as~$T\to+\infty$.

Accordingly, to complete the proof of~\eqref{UNST-2}, it only remains to
check that~$S_T\to1$ as~$T\to+\infty$.
To this end, we note that, by~\eqref{fanbdg} and~\eqref{Pop},
\begin{eqnarray*}
&&0=\frac1{2\pi c(S_T;T)\,h(S_T)^\frac{1-2S_T}{2S_T}}\,\frac{d}{ds}{\mathcal{E}}_2(S_T;T)\ge
\frac{f(S_T)\,h(S_T)}{c(S_T;T)}-S_T(2S_T-1)C_1
\\&&\qquad\quad=h(S_T)\Big((2S_T-1) \ln T
+(2S_T-1)\gamma-2S_T -2(2S_T-1)\Big)-S_T(2S_T-1)C_1\\&&\qquad\quad\ge
h(S_T)\frac{\ln T}2-C_2-C_3h(S_T),
\end{eqnarray*}
for some constants~$C_2$, $C_3>0$.

For this reason,
$$ \lim_{T\to+\infty} \frac{C_2}{h(S_T)}+C_3\ge \lim_{T\to+\infty}\frac{\ln T}{2}=+\infty,$$
whence
$$ 0=\lim_{T\to+\infty}h(S_T)=\lim_{T\to+\infty}
-2\zeta(-2s_T).$$
Since the only zero of the Riemann Zeta Function in~$[-2,-1]$
occurs at~$-2$, we thereby infer that~$-2s_T\to-2$,
and thus~$s_T\to1$, as~$T\to+\infty$.
With this, we have concluded the analytical verification of~\eqref{UNST-2}.

\subsection{Analytical verification of~\eqref{SLL-1}}\label{ALLc-1}

We compute that 
\begin{equation*}
\begin{split}
\frac{d}{ds}\mathcal{G}_1(s;L,T)&=\frac{T}{2\pi}\left[\frac{2\Gamma'(1+2s)\sin(\pi s)+\pi\Gamma(1+2s)\cos(\pi s)
-2\Gamma(1+2s)\sin(\pi s)\ln L}{L^{1+2s}}\right],\\
&=\frac{T\Gamma(1+2s)}{\pi L^{1+2s}}\left[\sin(\pi s)\big(\psi(1+2s)-\ln L\big)+\frac{\pi}{2}\cos(\pi s)\right].
\end{split}
\end{equation*}
For~$L$ large enough, we define 
\begin{equation*}
s_L^{(-)}:=\frac{1}{8 \ln L}
\end{equation*} 
and we notice that
\begin{equation}\label{c9and51.}
\lim_{L\to +\infty} s_L^{(-)}=0.
\end{equation} 
Also, if $L$ is sufficiently large, using \eqref{RViNSXC} (and noticing that~$Z(1+2s)<0$), we obtain that
\begin{equation*}
\begin{split}
\frac{\pi L^{1+2s}}{T\Gamma(1+2s)}\frac{d}{ds}\mathcal{G}_1(s;L,T)&> 
\sin(\pi s)\left(-\gamma-\frac{1}{1+2s}-\ln L\right)+\frac{\pi}{2}\cos(\pi s)\\
&\geq -2\sin(\pi s)\ln L+\frac{\pi}{2}\cos(\pi s)\\
&\geq -2\pi s\ln L+ \frac{\pi}{4}\\
&\geq 0,
\end{split}
\end{equation*}
for all $s\in \left(0,s_L^{(-)}\right]$. 

Analogously, we define 
\begin{equation}
s_L^{(+)}:=\frac{2}{3(\ln L-\psi(3))},
\end{equation}
and note that also in this case we have that 
\begin{equation}\label{c9and51.2}
\lim_{L\to +\infty} s_L^{(+)}=0.
\end{equation}
Now, we notice that if~$L$ is sufficiently large and~$s\in\left[\frac12,1\right)$,
$$ \frac{\pi L^{1+2s}}{T\Gamma(1+2s)}\frac{d}{ds}\mathcal{G}_1(s;L,T)\le-\frac{\sin(\pi s)}2 \ln L+\frac\pi2\cos(\pi s)<0.
$$
If instead~$s\in\left[s_L^{(+)},\frac12\right)$,
\begin{eqnarray*}
\frac{\pi L^{1+2s}}{T\Gamma(1+2s)}\frac{d}{ds}\mathcal{G}_1(s;L,T)&<&
\sin(\pi s)\big(\psi(3)-\ln L\big)+\frac{\pi}{2}\cos(\pi s)\\&\le&
-\sin(\pi s)\big(\ln L-\psi(3)\big)+\frac{\pi}{2}\\&\le&
-\sin\big(\pi s_L^{(+)}\big)\big(\ln L-\psi(3)\big)+\frac{\pi}{2}
\\&\le&-\frac{3\pi s_L^{(+)}}4\big(\ln L-\psi(3)\big)+\frac{\pi}{2}
\\&\le&0.
\end{eqnarray*}
In conclusion, if~$L$ is large enough and~$s\in\left[s_L^{(+)},1\right)$ we have that
\begin{equation*}
\frac{d}{ds}\mathcal{G}_1(s;L,T)<0.
\end{equation*}

Therefore, from these considerations we conclude that 
\begin{equation}\label{F3con-1.3}
\begin{split}&{\mbox{there exists at least one critical
point for~${\mathcal{G}}_1(s;L,T)$}}\\&{\mbox{and all the critical points of~${\mathcal{G}}_1(s;L,T)$
are located in~$\left(s^{(-)}_L,s^{(+)}_L\right)$.}}\end{split}\end{equation}

Now, we show that if $L$ is large enough, then   
\begin{equation}\label{F3con-1277}
{\mbox{$\displaystyle\frac{d}{ds}\mathcal{G}_1(s;L,T)$ is strictly decreasing in~$\left(s^{(-)}_L,s^{(+)}_L\right)$.}}
\end{equation}
To prove it, if $L$ is large enough, we observe that
\begin{equation*}
\begin{split}
&\frac{d^2}{d s^2}\mathcal{G}_1(s;L,T)\\
=\;&\frac{2T\Gamma(1+2s)}{\pi L^{1+2s}}\left(\sin(\pi s)\left(\psi(1+2s)-\ln L\right)+\frac{\pi}{2}\cos(\pi s)    \right)\left(\psi(1+2s)-\ln L\right)\\
&\qquad
+\frac{T\Gamma(1+2s)}{\pi L^{1+2s}}\left(\pi\cos(\pi s)\left(\psi(1+2s)-\ln L\right)+2\sin(\pi s)\psi'(1+2s)-\frac{\pi^2}{2}\sin(\pi s)\right)\\
=\;&\frac{2T\Gamma(1+2s)}{\pi L^{1+2s}}\Big(\sin(\pi s)\left(\psi(1+2s)-\ln L\right)+\pi\cos(\pi s)    \Big)\big(\psi(1+2s)-\ln L\big)\\
&\qquad +\frac{T\Gamma(1+2s)}{\pi L^{1+2s}}\left(2\sin(\pi s)\psi'(1+2s)-\frac{\pi^2}{2}\sin(\pi s)\right)\\
\leq\;& \frac{2T\Gamma(1+2s)}{L^{1+2s}}\Big( s\left(\psi(1+2s)-\ln L-\pi\right)+1\Big)\left(\psi(1+2s)-\ln L\right)\\
&\qquad+\frac{2 T\Gamma(1+2s)}{\pi L^{1+2s}}\sin(\pi s)\psi'(1+2s).
\end{split}
\end{equation*}

Now we observe that if~$s\in \left(s_L^{(-)},s_L^{(+)}\right)$ then
\begin{eqnarray*}
s\big(\ln L-\psi(1+2s)+\pi\big)-1&\le& \frac{2}{3(\ln L-\psi(3))}
\big(\ln L-\psi(1+2s)+\pi\big)-1\\&=&
\frac{2}{3}+\frac{2(\psi(3)-\psi(1+2s)+\pi)}{3(\ln L-\psi(3))}-1\\&=&-\frac13+\frac{2(\psi(3)-\psi(1+2s)+\pi)}{3(\ln L-\psi(3))}
\\&\le&-\frac16,
\end{eqnarray*}
as long as~$L$ is taken sufficiently large.

Moreover, since $\psi(1+2s)$ is increasing for $s\in (0,1)$, if~$L$ is sufficiently large,
we have that, for all~$s\in \left(s_L^{(-)},s_L^{(+)}\right)$,  
\begin{equation*}
\sin(\pi s)\psi'(1+2s)\leq \pi s\psi'(1+2s)\leq \pi s_L^{(+)}\psi'(1+2s).
\end{equation*}
Thus, using these pieces of information, we conclude that, if $L$ is large enough, for every~$s\in \left(s_L^{(-)},s_L^{(+)}\right)$,
\begin{eqnarray*}
\frac{d^2}{d s^2}\mathcal{G}_1(s;L,T)&\leq& -
\frac{T\Gamma(1+2s)}{3L^{1+2s}}\left(\ln L-\psi(1+2s)\right)+\frac{2 T\Gamma(1+2s)}{ L^{1+2s}} s_L^{(+)}\psi'(1+2s)\\
&=& -
\frac{T\Gamma(1+2s)}{3L^{1+2s}}\big(\ln L-\psi(1+2s)-6s_L^{(+)}\psi'(1+2s)\big)
\\&=&-\frac{T\Gamma(1+2s)}{3L^{1+2s}}\left(\ln L-\psi(1+2s)-\frac{4\psi'(1+2s)}{\ln L-\psi(3)}\right)
\\&<&0.\end{eqnarray*}
This establishes the claim in~\eqref{F3con-1277}.

Now, the claim in~\eqref{SLL-1} follows from \eqref{F3con-1.3}, \eqref{F3con-1277} and the limits in~\eqref{c9and51.} and~\eqref{c9and51.2}.

\subsection{Analytical verification of \eqref{SLL-1.K}}\label{ALLc-1.kk}

We compute that 
\begin{equation*}
\begin{split}
\frac{d}{ds}\mathcal{G}_2(s;L,T)&=-\frac{T}{4}\left[\frac{2\ln L}{L^{1+2s}\zeta(1+2s)}
+\frac{2\zeta'(1+2s)}{L^{1+2s}\zeta^2(1+2s)}\right]\\
&=-\frac{T}{2L^{1+2s}\zeta(1+2s)}\left[\ln L+\frac{\zeta'(1+2s)}{\zeta(1+2s)}\right].
\end{split}
\end{equation*}
The Riemann zeta function $\zeta$ is a meromorphic function with pole at~$z=1$. We also recall that near~$z=1$ the derivative of the Riemann Zeta Function has the Laurent expansion
\begin{equation} \label{LAUR}\frac{\zeta'(z)}{\zeta(z)}=
-\frac{1}{z-1}+\gamma+O(z-1),\end{equation}
see e.g. page~481 in~\cite{CHO}. This leads to the expression 
\begin{equation}\label{fa-rr-rr}
\frac{\zeta'(1+2s)}{\zeta(1+2s)}=-\frac{1}{2s}+\gamma+O(2s)
\end{equation}
and to the existence of some constant~$C\in \R$ such that 
\begin{equation}\label{s3icrub485yc45yuxj}
\frac{\zeta'(1+2s)}{\zeta(1+2s)}\geq -\frac{1}{2s}+\gamma+ C s,
\end{equation}
for all $s\in (0,1)$. We point out that we can assume without loss of generality that $C\in (-\infty,0]$.

Therefore, if we define  
\begin{equation*}
s_L^{(+)}:= \frac{1}{2\left(\ln L+\gamma+C\right)},
\end{equation*}
we obtain that 
\begin{equation}\label{fitsikhif}
\begin{split}
\frac{2 L^{1+2s}\zeta(1+2s)}{T}\frac{d}{ds}\mathcal{G}_2(s;L,T)&=-\left[\ln L+\frac{\zeta'(1+2s)}{\zeta(1+2s)}\right]\\
&\leq -\left(\ln L-\frac{1}{2s}+\gamma+ C s  \right)\\
&< -\left(\ln L-\frac{1}{2s}+\gamma+ C  \right)\\
&\le 0,
\end{split}
\end{equation} 
for all $s\in \left[s_{L}^{(+)},1\right)$.

We notice also that 
\begin{equation}\label{c9and51.cc}
\lim_{L\to +\infty} s_{L}^{(+)}=0.
\end{equation} 

Analogously, thanks to the expansion in \eqref{fa-rr-rr}, we obtain the existence of some~$\tilde{C}\in \R$ such that 
\begin{equation}\label{bi-prltr}
\frac{\zeta'(1+2s)}{\zeta(1+2s)}\leq -\frac{1}{2s}+\gamma+ \tilde{C}s,
\end{equation} 
for all $s\in (0,1)$. Without loss of generality, we can suppose that~$\tilde{C}\in [0,+\infty)$.

Thus, if we define  
\begin{equation*}
s_L^{(-)}:=\frac{1}{2(\ln L+\gamma+\tilde{C})},
\end{equation*}
we deduce that
\begin{equation}\label{tsndtsnd}
\begin{split}
\frac{2 L^{1+2s}\zeta(1+2s)}{T}\frac{d}{ds}\mathcal{G}_2(s;L,T)&=-\left[\ln L+\frac{\zeta'(1+2s)}{\zeta(1+2s)}\right]\\
&\geq -\left(\ln L-\frac{1}{2s}+\gamma+ \tilde{C} s \right)\\
&> -\left(\ln L-\frac{1}{2s}+\gamma+ \tilde{C} \right)\\
&\ge 0,
\end{split}
\end{equation}
for all $s\in \left(0,s_L^{(-)}\right]$.

We also observe that 
\begin{equation}\label{c9and51.2.cc}
\lim_{L\to +\infty} s_L^{(-)}=0.
\end{equation}

From \eqref{fitsikhif} and \eqref{tsndtsnd} we deduce that
\begin{equation}\label{F3con-1.3.cc}
\begin{split}
&{\mbox{there exists at least one critical
point for~${\mathcal{G}}_2(s;L,T)$}}\\&{\mbox{and all the critical points of~${\mathcal{G}}_2(s;L,T)$ are located in~$\left(s^{(-)}_L,s^{(+)}_L\right)$.}}
\end{split}
\end{equation}
We now show that, for all $L$ large enough,
\begin{equation}\label{F3con-1277.cc}
{\mbox{$\displaystyle{\frac{d}{ds}\mathcal{G}_2}(s;L,T)$ is strictly decreasing in~$\left(s^{(-)}_L,s^{(+)}_L\right)$.}}
\end{equation}
To prove it, using \eqref{bi-prltr} we estimate the derivative 
\begin{equation}\label{sjeiwoKKKKKKKruoh98o098765400}
\begin{split}
\frac{d^2}{ds^2}\mathcal{G}_2(s;L,T)&=\frac{T}{\zeta(1+2s)L^{1+2s}}\left[\ln L+\frac{\zeta'(1+2s)}{\zeta(1+2s)}\right]^2\\
&\qquad-\frac{T}{\zeta(1+2s)L^{1+2s}}\left[\frac{\zeta''(1+2s)}{\zeta(1+2s)}-\frac{(\zeta'(1+2s))^2}{\zeta^2(1+2s)}\right].
\end{split}
\end{equation}
Now we notice that
\begin{equation}\label{sjwieru3HHHH4ity5iuik896okyhtg}
\left| \ln L+\frac{\zeta'(1+2s)}{\zeta(1+2s)}
\right| \le \ln L -\frac1{2s}+\gamma +\max\{\tilde C, -C\}s,\end{equation}
thanks to~\eqref{s3icrub485yc45yuxj} and~\eqref{bi-prltr}.

Moreover, by the definitions of~$s_L^{(-)}$ and~$s_L^{(+)}$, we have that
$$ \frac{1}{2(\ln L+\gamma+\tilde{C})}<s<\frac{1}{2\left(\ln L+\gamma+C\right)},$$
which gives that
$$ \left| \ln L-\frac1{2s}+\gamma \right| \le \max\{\tilde C, -C\}.$$
{F}rom this and~\eqref{sjwieru3HHHH4ity5iuik896okyhtg}, we deduce that
\begin{equation}\label{sjeiwoKKKKKKKruoh98o0987654}\left| \ln L+\frac{\zeta'(1+2s)}{\zeta(1+2s)}
\right| \le \max\{\tilde C, -C\}(1+s).
\end{equation}

Additionally, differentiating the Laurent expansion in~\eqref{LAUR},
\begin{eqnarray}\label{sl-cd.99}
\frac{\zeta''(z)}{\zeta(z)}-
\left(\frac{\zeta'(z)}{\zeta(z)}\right)^2=
\frac{d}{dz}\frac{\zeta'(z)}{\zeta(z)}=\frac{1}{(z-1)^2}+O(1),
\end{eqnarray}
which leads to the existence of some constant $\hat{C}\in \mathbb{R}$ such that  
\begin{equation}\label{sjeiwoKKKKKKKruoh98o0987654BIS}
\frac{\zeta''(1+2s)}{\zeta(1+2s)}-
\left(\frac{\zeta'(1+2s)}{\zeta(1+2s)}\right)^2 \geq \frac{1}{(2s)^2}+\hat{C},
\end{equation}
for all $s\in (0,1)$. 

Using~\eqref{sjeiwoKKKKKKKruoh98o0987654} and~\eqref{sjeiwoKKKKKKKruoh98o0987654BIS}
into~\eqref{sjeiwoKKKKKKKruoh98o098765400}, we conclude that
\begin{eqnarray*}&&
\frac{d^2}{ds^2}\mathcal{G}_2(s;L,T)\\&\leq&
\frac{T}{\zeta(1+2s)L^{1+2s}}\left[(\max\{\tilde C, -C\})^2(1+s)^2- \frac{1}{(2s)^2}-\hat{C}\right]
\\&\le&
\frac{T}{\zeta(1+2s)L^{1+2s}}\left[(\max\{\tilde C, -C\})^2 \left(1+\frac{1}{2\left(\ln L+\gamma+C\right)}\right)^2
-\big(\ln L+\gamma+C\big)^2-\hat{C}\right],
\end{eqnarray*}
for all~$s\in \left(s_L^{(-)},s_L^{(+)}\right)$.

Thus, if~$L$ is sufficiently large, we obtain~\eqref{F3con-1277.cc}.

Therefore, claim \eqref{SLL-1.K} follows from \eqref{F3con-1.3.cc}, \eqref{F3con-1277.cc} and the limits in \eqref{c9and51.cc} and \eqref{c9and51.2.cc}.

\subsection{Analytical verification of~\eqref{SLL}}\label{ALL}

We observe that
\begin{eqnarray*}&&
{\mathcal{P}}(s;L,T)\\&:=&
\frac{{\pi}\,L^{1+2s}\, \zeta(2s) }{T\Gamma(1+2s)}\;
\frac{d}{ds}{\mathcal{G}}_3(s;L,T)
\\&=&
\zeta(1+2 s) \Big(
{- \ln L} \,\sin(\pi s) + {\frac{\pi}2} \cos(\pi s)
+ {\sin(\pi s)} \psi(1+2 s )\Big) +  {\sin(\pi s)} \zeta'(1+2 s)\\
&&\qquad- \frac{{\zeta(1+2 s)} \sin(\pi s) \zeta'(2 s)}{\zeta(2s)}\end{eqnarray*}
and the positivity of the derivative of~${\mathcal{G}}_3$ is equivalent to the positivity of~${\mathcal{P}}$.

Let also
$$ {\mathcal{Q}}(s):=
\zeta(1+2 s) \left(
{\frac{\pi}2} \cos(\pi s) + {\sin(\pi s)} \psi(1+2 s)\right) + {\sin(\pi s)} \zeta'(1+2 s)$$
and note that~$\sup_{s\in\left(\frac12,1\right)}|{\mathcal{Q}}(s)|\le C$ for some~$C>0$;
in addition,
\begin{eqnarray*}
{\mathcal{P}}(s;L,T)&=&
{-\zeta(1+2 s)} \ln L \,\sin(\pi s) - \frac{{\zeta(1+2 s)} \sin(\pi s) \zeta'(2 s)}{\zeta(2s)}+
{\mathcal{Q}}(s)\\
&=&{- \zeta(1+2 s)} \sin(\pi s)\left(
\ln L+ \frac{ \zeta'(2 s)}{\zeta(2s)}\right)+
{\mathcal{Q}}(s)\\&=&{- \zeta(1+2 s)} \sin(\pi s)\left(
\ln L-\frac{1}{2s-1}\right)+
\widetilde{\mathcal{Q}}(s),
\end{eqnarray*}
where
$$ \widetilde{\mathcal{Q}}(s):={\mathcal{Q}}(s)
{- \zeta(1+2 s)} \sin(\pi s)\left(
\frac{ \zeta'(2 s)}{\zeta(2s)}+\frac{1}{2s-1}\right).$$
We stress that~$\sup_{s\in\left(\frac12,1\right)}|\widetilde {\mathcal{Q}}(s)|\le \widetilde C$ for some~$\widetilde C>0$, thanks to~\eqref{LAUR}.

Let now
$$ \varepsilon_L:=\frac{1}{\sqrt{\ln L}}\qquad{\mbox{ and }}\qquad
s^{(\pm)}_L:=\frac12+\frac{1\pm \varepsilon_L}{2\ln L}.$$
We notice that
\begin{equation}\label{SELE} \lim_{L\to+\infty}\varepsilon_L=0
\qquad{\mbox{and}}\qquad\lim_{L\to+\infty}s^{(\pm)}_L=\frac12.\end{equation}
Also, when~$s\in\left(\frac12,s^{(-)}_L\right]$, we have that
$$ 2s-1\le 2s^{(-)}_L-1=
\frac{1- \varepsilon_L}{\ln L}$$
and, as a result, for large~$L$,
\begin{eqnarray*}
{\mathcal{P}}(s;L,T)&\ge&{\zeta(1+2 s)} \sin(\pi s)\left(
-\ln L+\frac{\ln L}{1- \varepsilon_L}
\right)-
\widetilde{C}
\\&=&\frac{{ \varepsilon_L}\,\ln L\;\zeta(1+2 s) \sin(\pi s)}{1- \varepsilon_L}-
\widetilde{C}
\\&\geq&\frac{ \sqrt{\ln L}\;\zeta(1+2 s) \sin(\pi s)}{{2}}-
\widetilde{C}\\&\geq& \frac{\sqrt{\ln L}\;\zeta(3)}{{4}} -
\widetilde{C}\\&>&0.
\end{eqnarray*}

Instead, when~$s\in\left[s^{(+)}_L,1\right)$,
$$ 2s-1\ge 2s^{(+)}_L-1=
\frac{1+ \varepsilon_L}{\ln L},$$
which entails that
\begin{eqnarray*}
{\mathcal{P}}(s;L,T)&\le&{ \zeta(1+2 s)} \sin(\pi s)\left(
-\ln L+\frac{\ln L}{1+ \varepsilon_L}\right)+
\widetilde{\mathcal{Q}}(s)\\&\le&
-\frac{{\varepsilon_L}\,\ln L\;\zeta(1+2 s) \sin(\pi s)}{1+ \varepsilon_L}
+{\mathcal{Q}}(s)
+C_1 \sin(\pi s)\\&\le&
-\frac{{ \sqrt{\ln L}}\;\zeta(3) \sin(\pi s)}{1+ \varepsilon_L}
+\frac{\pi}{{2}} \cos(\pi s)\zeta(1+2 s) 
+C_2 \sin(\pi s)\\
&=&
-\sin(\pi s)\left(
\frac{{\sqrt{\ln L}}\;\zeta(3) }{1+ \varepsilon_L}-C_2\right)
+\frac{\pi}{{2}} \cos(\pi s)\zeta(1+2 s) \\
&\le&
-\sin(\pi s)\left(\frac{ \sqrt{\ln L}\;\zeta(3)}{{2}} -C_2\right)
+\frac{\pi}{{2}} \cos(\pi s)\zeta(1+2 s) 
\\&<&0,
\end{eqnarray*}
where~$C_1$ and~$C_2$
are suitable positive constants.

{F}rom these considerations, it follows that
there exists at least one zero for~${\mathcal{P}}$
and all the zeros of~${\mathcal{P}}$ are located in~$\left(s^{(-)}_L,s^{(+)}_L\right)$.
As a consequence,
\begin{equation}\label{F3con-1}
\begin{split}&{\mbox{there exists at least a critical
point for~${\mathcal{G}}_3$}}\\&{\mbox{and all the critical points of~${\mathcal{G}}_3$ are located in~$\left(s^{(-)}_L,s^{(+)}_L\right)$.}}\end{split}\end{equation}

We now show that
\begin{equation}\label{F3con}
{\mbox{${\mathcal{P}}$ is strictly decreasing in~$\left(s^{(-)}_L,s^{(+)}_L\right)$.}}
\end{equation}
For this, we calculate that
\begin{equation}\label{SEMDsiakd}
\begin{split}
\frac{d}{ds}{\mathcal{P}}(s;L,T)
\,&=\,
\ln L\, \Big(-{2} \sin(\pi s) \zeta'(1+2 s) - {\pi} \zeta(1+2 s) \cos(\pi s)\Big)
\\&\quad+ \pi \cos(\pi s) \left[{2}\zeta'(1+2 s) + {\zeta(1+2 s)} 
\left( \psi(1+2 s ) - \frac{ \zeta'(2 s)}{\zeta(2 s)}\right)\right]
\\&\quad+ {2}\sin(\pi s) \left[ \zeta''(1+2 s ) + \zeta'(1+2 s) \left( \psi(1+2 s) 
-\frac{ \zeta'(2 s)}{\zeta(2 s)}\right) \right.\\&\quad\qquad\left.
+ \zeta(1+2 s) \left(-\frac{ \zeta''(2 s)}{
\zeta(2 s)} 
+ \left(\frac{\zeta'(2 s)}{\zeta(2 s)}\right)^2 + \psi'(1+2 s) - \pi^2\right)\right]\\
&\le\,
\ln L\, \Big(-{2} \sin(\pi s) \zeta'(1+2 s) - {\pi }\zeta(1+2 s) \cos(\pi s)\Big)
\\&\quad- \pi \cos(\pi s) \frac{ \zeta'(2 s)}{\zeta(2 s)}
- {2}\sin(\pi s) \frac{ \zeta'(2 s)}{\zeta(2 s)} \\&\quad
+ \zeta(1+2 s) \left(-\frac{ \zeta''(2 s)}{
\zeta(2 s)} 
+ \left(\frac{\zeta'(2 s)}{\zeta(2 s)}\right)^2 \right)+C_3
\end{split}\end{equation}
for some constant~$C_3>0$.

Plugging~\eqref{sl-cd.99} into~\eqref{SEMDsiakd} we find that
\begin{equation*}
\begin{split}
\frac{d}{ds}{\mathcal{P}}(s;L,T)
\,&\le\,
\ln L\, \Big(-{2} \sin(\pi s) \zeta'(1+2 s) - { \pi} \zeta(1+2 s) \cos(\pi s)\Big)
\\&\quad- \pi \cos(\pi s) \frac{ \zeta'(2 s)}{\zeta(2 s)}
- {2}\sin(\pi s) \frac{ \zeta'(2 s)}{\zeta(2 s)}
-\frac{ \zeta(1+2 s) }{(2s-1)^2}+C_4
\end{split}\end{equation*}
for some constant~$C_4>0$. Thus, since~$\zeta' \le0$,
\begin{equation}\label{LSUNIBS}
\begin{split}
\frac{d}{ds}{\mathcal{P}}(s;L,T)
\,&\le\,
\ln L\, \Big(-{2}\sin(\pi s) \zeta'(1+2 s) - {\pi} \zeta(1+2 s) \cos(\pi s)\Big)
\\&\quad\qquad
- \frac{ {2}\zeta'(2 s)}{\zeta(2 s)}
-\frac{ \zeta(1+2 s) }{(2s-1)^2}+C_4\\&\le\,
\ln L\, \Big(-{2} \sin(\pi s) \zeta'(1+2 s) - {\pi} \zeta(1+2 s) \cos(\pi s)\Big)
\\&\quad\qquad
+\frac{{2}}{2s-1}
-\frac{ \zeta(1+2 s) }{(2s-1)^2}+C_5
\end{split}\end{equation}
for some~$C_5>0$, where~\eqref{LAUR} has been used once again.

Moreover, if~$s\in\left(\frac12,s^{(+)}_L\right]$, then, for large~$L$,
$$ 2s-1\le
\frac{1+\varepsilon_L}{\ln L}\le\frac{2}{\ln L}$$
and therefore we deduce from~\eqref{LSUNIBS} that
\begin{eqnarray*}&&
\frac{d}{ds}{\mathcal{P}}(s;L,T)\leq
\frac{C_6}{2s-1}
-\frac{ \zeta(1+2 s) }{(2s-1)^2}\leq
\frac{1}{2s-1}\left(
C_6-\frac{ \zeta(3) }{2s-1}
\right)\\&&\qquad\qquad\leq-
\frac{1}{2s-1}\left(
\frac{ \zeta(3) }{2s^{(+)}_L-1}-C_6
\right)\leq-
\frac{1}{2s-1}\left(
\frac{ \zeta(3) \,\ln L}{2}-C_6
\right)<0,
\end{eqnarray*}
which completes the proof of~\eqref{F3con}.

The desired claim in~\eqref{SLL} now follows
combining~\eqref{F3con-1}
and~\eqref{F3con}
with the second limit in~\eqref{SELE}.

\subsection{Analytical verification of~\eqref{ASOG4}}\label{A4444}

We start by computing the derivative of $\mathcal{G}_4$ with respect to $s$, and we get that
\begin{equation}\label{DerivG4}
\frac{d}{ds}\mathcal{G}_4(s,T,L)=-\frac{T}{{2}L^{1+2s}\,\zeta(2s)}\left(\ln L+\frac{\zeta'(2s)}{\zeta(2s)}\right).
\end{equation}
Now we observe that, for $s\in\left(\frac{1}{2},1\right)$, the function
$$
m(s):=\frac{\zeta'(2s)}{\zeta(2s)}
$$
is negative and strictly increasing. Furthermore, from \eqref{LAUR} we infer that 
\begin{equation}\label{SPSPSPS}
\lim_{s\searrow1/2}m(s)=-\infty.
\end{equation}
Now, in light of~\eqref{DerivG4},
\begin{equation}\label{KM:LD 03l4fnnbvkf8586}
\begin{split} &\frac{d}{ds}\mathcal{G}_4(s,T,L)>0{\mbox{ if and only if~$\ln L+m(s)<0$,}}\\&
{\mbox{and thus if and only if~$L < \exp\left( -m(s)\right)$.}}\end{split}\end{equation}
Also, from~\eqref{ELLSTARVO} and the monotonicity of the function~$m$,
for each~$L\le L^\star$ and~$s\in\left(\frac12,1\right)$ we have that
$$ \exp\left( -m(s)\right)>\exp\left( -m(1)\right)=L^\star\ge L,$$
and therefore we deduce from~\eqref{KM:LD 03l4fnnbvkf8586} that
\begin{equation}\label{REVCD-BAN3JNDMN}
{\mbox{when~$L\leq L^\star$ the supremum of~$\left(\frac12,1\right)\ni
s\mapsto {\mathcal{G}}_4(s;L,T)$ is uniquely attained at~$s=1$.}}
\end{equation}
If instead~$L>L^\star$, using~\eqref{SPSPSPS}
we see that
\begin{eqnarray*}
&&\lim_{s\searrow1/2}\exp\left( -m(s)\right)=+\infty>L\\
\\{\mbox{and }}&&\lim_{s\nearrow1}\exp\left( -m(s)\right)=
\exp\left( -m(1)\right)=L^\star<L.\end{eqnarray*}
This and the strict monotonicity of~$m$
yield that for each $L>L^\star$ there exists a unique $s_L\in\left(\frac{1}{2},1\right)$ such that 
\begin{equation}\label{ZerodidG}
\frac{d}{ds}\mathcal{G}_4(s_L,L,T)=0,
\end{equation}
with~$\frac{d}{ds}\mathcal{G}_4(s,L,T) > 0$ if and only if~$  s \in \left(\frac{1}{2},s_L\right)$,
namely~$s_L$ is the unique maximum for $\mathcal{G}_4(s,L,T)$ when $L>L^\star$. 

{F}rom these observations and~\eqref{REVCD-BAN3JNDMN},
in order to complete the analytical proof of \eqref{ASOG4}, it is only left to show that
\begin{equation}\label{FALRE}
\lim_{L\rightarrow+\infty}s_L=\frac{1}{2}.
\end{equation} To this end, equations \eqref{DerivG4} and \eqref{ZerodidG} give that
\begin{equation*}
0=\ln L+m(s_L),
\end{equation*}
which leads to
\begin{equation}
+\infty=\lim_{L\to+\infty}\ln L=-\lim_{L\to+\infty}m(s_L).
\end{equation}
Since the only pole of $m(s)$ at $\left[\frac{1}{2},1\right]$
occurs in $s=\frac{1}{2}$, then we obtain~\eqref{FALRE}, as desired.

\end{appendix}

\section*{Conclusions}
In this article we have introduced several efficiency functionals accounting for different foraging strategies of a predator (for simplicity, in a one-dimensional environment). The foraging strategies rely on the possible adoption of a different L\'evy exponent for the diffusion of the predator. The corresponding efficiency functionals thus compare, roughly speaking, the outcome of the forager's hunt with the effort required to implement it and our objective was to understand, in different scenarios, which L\'evy exponent optimizes, or pessimizes, a given efficiency functional for a given distribution of targets.

Several biological environments have been considered, such as the ones of a target close to the initial position of the predator, of a sparse distribution of targets, of a remotely located target, etc.

The optimal exponents correspond to different hunting strategies, driven by either classical or anomalous types of diffusion: specifically, one can compare the classical diffusion induced by Brownian motion and the one obtained by an inverse square law distribution. Our findings show that the {\em theoretical optimality} of an exponent has to be confronted with the {\em practical reliability} of the corresponding hunting strategy, since we pointed out the occurrence of {\em bifurcation phenomena} depending on the environmental parameters in which a sudden switch takes place between a theoretically optimal exponent and a less ideal, but more secure, foraging procedure. In particular, optimal exponents can be located arbitrarily close to pessimal ones, making the practical choice of the hunting strategy a delicate balance between a striving for maximum success and a more conservative attitude to prioritize safety.

Cases of intermediate optimal exponents, balancing between purely Gaussian and L\'evy distributions, have also been detected.

We kept our analysis as simple as possible, by removing additional parameters in the model such as the possibility of restarting the strategy after each hit and the distinction between random relocation and direct finite-range vision. Also, we did not introduce any additional a-priori bound on the length of the forager's journey and no additional truncation on the power law distributions has been added into the model. Furthermore, the efficiency functionals did not account just for a single foraging success but rather for its time average.

In our computation, a significant role is also played by the possibly different normalization constants involved in probabilistic and analytical models. The results obtained present explicit solutions in closed form which utilize only elementary special functions (in particular, no expensive numerical simulations were needed and the results are not affected by truncations or discretizations).

\end{document}